\documentclass[a4paper,10pt]{article}
\usepackage[top=13mm, bottom=23mm, inner=20mm, outer=20mm]{geometry}

\usepackage{mathtools,amssymb,mathrsfs,bm}

\usepackage{graphicx}
\usepackage{caption,subcaption,float}
\usepackage{tabularx,multirow}
\usepackage[inline]{enumitem}

\usepackage{xcolor}

\numberwithin{equation}{section}
\numberwithin{figure}{section}
\numberwithin{table}{section}

\usepackage{calc,etoolbox}

\usepackage[inline]{enumitem}

\makeatletter
\@ifpackageloaded{hyperref}{}{\usepackage[hidelinks]{hyperref}}
\@ifpackageloaded{cleveref}{}{
	\usepackage[compress]{cleveref}

	\crefformat{section}{\S#2#1#3}
	\crefformat{subsection}{\S#2#1#3}
	\crefformat{subsubsection}{\S#2#1#3}
	\crefrangeformat{section}{\S\S#3#1#4 to~#5#2#6}
	\crefmultiformat{section}{\S\S#2#1#3}{ and~#2#1#3}{, #2#1#3}{ and~#2#1#3}
	
	\crefname{figure}{fig.}{figs.}
	\Crefname{figure}{Fig.}{Figs.}
	
	\crefname{equation}{}{}
	\crefformat{equation}{(#2#1#3)}
	\crefmultiformat{equation}{(#2#1#3)}{ and (#2#1#3)}{, (#2#1#3)}{ and (#2#1#3)}
}
\makeatother

\usepackage[title]{appendix}

\crefname{appsec}{Appendix}{Appendices}
\Crefname{appsec}{Appendix}{Appendices}

\usepackage{xspace}
\makeatletter
\DeclareRobustCommand\onedot{\futurelet\@let@token\@onedot}
\def\@onedot{\ifx\@let@token.\else.\null\fi\xspace}
\makeatother
\def\eg{e.g\onedot} 
\def\ie{i.e\onedot}

\def\etal{et al\onedot}
\def\sf{s.f\onedot}

\makeatletter
\newcommand{\floattag}[1]{%
  \@namedef{the\@captype}{#1}%
  \@namedef{theH\@captype}{#1}%
  \addtocounter{\@captype}{-1}}
\makeatother
\input{Templates/StandardMathsMacros}
\usepackage{amsthm, thmtools, thm-restate}
\theoremstyle{plain}

\theoremstyle{definition}

\theoremstyle{remark}

\theoremstyle{definition}
\declaretheorem[name=Test problem, numberwithin=section]{testproblem}

\crefformat{testproblem}{test problem #2#1#3}
\crefrangeformat{testproblem}{test problems #3#1#4 to~#5#2#6}
\Crefrangeformat{testproblem}{Test problems #3#1#4 to~#5#2#6}
\crefmultiformat{testproblem}{test problems #2#1#3}{ and~#2#1#3}{, #2#1#3}{ and~#2#1#3}

\usepackage{doi}
\graphicspath{{Figures/}}

\begin{document}

\title{A Well Balanced Reconstruction with Bounded Velocities and Low-Oscillation Slow Shocks for the Shallow Water Equations}
\author{Edward W. G. Skevington}
\date{21/06/2021}
\maketitle

\begin{abstract}

Many numerical schemes for hyperbolic systems require a piecewise polynomial reconstruction of the cell averaged values, and to simulate perturbed steady states accurately we require a so called ‘well balanced’ reconstruction scheme. For the shallow water system this involves reconstructing in surface elevation, to which modifications must be made as the fluid depth becomes small to ensure positivity.

We investigate the scheme proposed in [Skevington] though numerical experiments, demonstrating its ability to resolve steady and near steady states at high accuracy. We also present a modification to the scheme which enables the resolution of slowly moving shocks and dam break problems without compromising the well balanced property.
\end{abstract}

%
%
%
%
%
%
%
%
%
%


\section{Introduction}

An important system of equations in fluid dynamical research is the shallow water equations, which enforce conservation of mass and momentum for a free surface flow under hydrostatic pressure. Including the effects topography and assuming the flow to be unidirectional results in the system (\eg \cite{bk_Stoker_WW,bk_Ungarish_GCI})
\begin{subequations}	\label{eqn:INTRO_SW_Full}	\begin{align}
	\pdv{h}{t} + \pdv{}{x} (uh) &= 0,		\label{eqn:INTRO_SW_Full_h}	\\
	\pdv{}{t}(uh) + \pdv{}{x} \ppar*{ u^2 h + \frac{g h^2}{2} } &= - g h\pdv{b}{x},			\label{eqn:INTRO_SW_Full_u}
\end{align}\end{subequations}
where $h(x,t)$ is the depth of the fluid, $b(x)$ the bed elevation, $u(x,t)$ the horizontal velocity, and $g$ the gravitational force per unit mass.

\begin{figure}[tp]
	\centering
\begingroup%
  \makeatletter%
  \providecommand\color[2][]{%
    \errmessage{(Inkscape) Color is used for the text in Inkscape, but the package 'color.sty' is not loaded}%
    \renewcommand\color[2][]{}%
  }%
  \providecommand\transparent[1]{%
    \errmessage{(Inkscape) Transparency is used (non-zero) for the text in Inkscape, but the package 'transparent.sty' is not loaded}%
    \renewcommand\transparent[1]{}%
  }%
  \providecommand\rotatebox[2]{#2}%
  \newcommand*\fsize{\dimexpr\f@size pt\relax}%
  \newcommand*\lineheight[1]{\fontsize{\fsize}{#1\fsize}\selectfont}%
  \ifx\svgwidth\undefined%
    \setlength{\unitlength}{325.98425197bp}%
    \ifx\svgscale\undefined%
      \relax%
    \else%
      \setlength{\unitlength}{\unitlength * \real{\svgscale}}%
    \fi%
  \else%
    \setlength{\unitlength}{\svgwidth}%
  \fi%
  \global\let\svgwidth\undefined%
  \global\let\svgscale\undefined%
  \makeatother%
  \begin{picture}(1,0.26086957)%
    \lineheight{1}%
    \setlength\tabcolsep{0pt}%
    \put(0,0){\includegraphics[width=\unitlength,page=1]{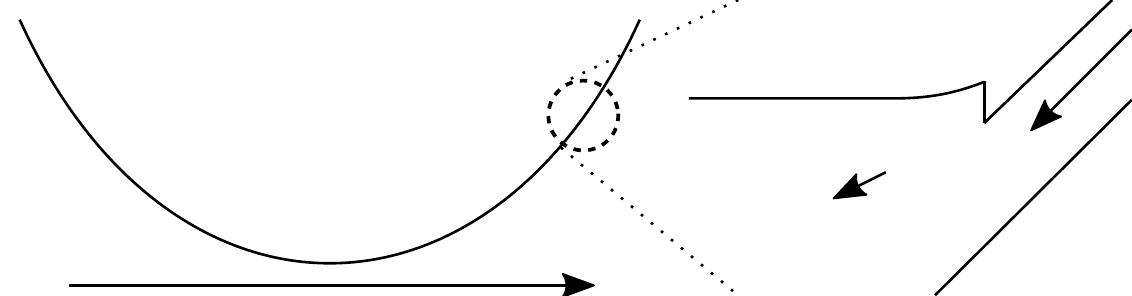}}%
    \put(0.51024246,0.01784084){\color[rgb]{0,0,0}\makebox(0,0)[lt]{\lineheight{1.25}\smash{\begin{tabular}[t]{l}$x$\end{tabular}}}}%
    \put(0.38542616,0.05262345){\color[rgb]{0,0,0}\makebox(0,0)[rt]{\lineheight{1.25}\smash{\begin{tabular}[t]{r}$b$\end{tabular}}}}%
    \put(0.28741057,0.18305797){\color[rgb]{0,0,0}\makebox(0,0)[t]{\lineheight{1.25}\smash{\begin{tabular}[t]{c}$\eta$\end{tabular}}}}%
    \put(0.76920047,0.11101976){\color[rgb]{0,0,0}\makebox(0,0)[rt]{\lineheight{1.25}\smash{\begin{tabular}[t]{r}$u$\end{tabular}}}}%
    \put(0,0){\includegraphics[width=\unitlength,page=2]{DrainIntoLake.pdf}}%
  \end{picture}%
\endgroup%

	\caption{Illustration of a typical situation, where the flow is approaching steady state leaving some fluid still on the slopes. The body of the fluid is close to `lake at rest', and on the slopes the fluid is close to `thin film'.}
	\label{fig:PWB_draining_flow}
\end{figure}

Numerical schemes which can accurately resolve steady states (or perturbations to them) must have a small error in the balance between the flux and source terms in the governing system, and schemes for which this balance is resolved at higher order than the formal order of the scheme are termed \emph{well balanced} \cite{ar_Bermudez_1994,ar_Greenberg_1996}. For \cref{eqn:INTRO_SW_Full} a special steady state is $h=0$, the \emph{dry bed}, whilst all other steady states form a family called the \emph{moving water equilibrium} for which the volume flux per unit width, $q$, and the energy per unit mass, $E$, are constant, where
\begin{align}	\label{eqn:SW_steady_moving}
	q &\eqdef uh,
	&
	E &\eqdef \frac{u^2}{2} + g \eta.
\end{align}
and $\eta \eqdef h+b$ is the surface elevation. There are two special cases of \cref{eqn:SW_steady_moving}. Firstly, the \emph{lake at rest} state (\cref{fig:PWB_draining_flow}) has $u = 0$ and $\eta$ constant. Secondly, the \emph{thin film} state (\cref{fig:PWB_draining_flow}) where the fluid is draining off the slope at a very high Froude number, \ie $\Fro \gg 1$ where $\Fro \eqdef \abs{u}/\sqrt{gh}$. In this state $h \ll H$ ($H$ the characteristic depth-scale) while $u = \Torder{\sqrt{gH}}$ ($\Torder{}$ is asymptotic notation), thus from $\dv*{q}{x} = \dv*{E}{x} = 0$ we deduce
\begin{align}
	\dv{h}{x} = \frac{1}{\Fro^2 - 1} \dv{b}{x} = \Torder{\frac{1}{\Fro^2} \frac{H}{L}} = \Torder{\frac{h}{L}}
\end{align}
($L$ the characteristic length-scale) and the variation of depth is small, whilst the variation in surface elevation is large, $\dv*{\eta}{x} = \Torder{H/L}$.

Well balanced schemes have been designed for a variety of numerical methods, \eg discontinuous Galerkin schemes \cite{ar_Noelle_2009,ar_Xing_2010,ar_Bollermann_2011,ar_Xing_2014,ar_Chen_2019} and the finite volume schemes we discuss here. These are well balanced on the basis of the reconstruction of the simulated fields \cite{ar_Kurganov_2002,ar_Kurganov_2007,ar_Bollermann_2013,ar_Chertock_2015,ar_Cheng_2016,ar_Cheng_2019}, and employing a well balanced numerical source \cite{ar_Bermudez_1994,ar_Bernstein_2016}. Recently, a framework for constructing well balanced schemes was proposed in \cite{ar_Skevington_N002_Recon_Theory}, where two reconstructions were performed, one in depth and the other in surface elevation. By performing a convex combination of the two reconstructions it is possible to enforce that the reconstructed depth in each cell increases monotonically with the cell average depth, thereby obtaining a lower bound on the reconstructed depth and an upper bound on the reconstructed velocity. Here we explore the scheme derived in \cite{ar_Skevington_N002_Recon_Theory} using these principles, and compare the resulting simulations to other well balanced reconstructions from the literature. The primary goal of this paper is to demonstrate the accuracy that can be achieved by the convex combination scheme.

We also tackle other challenges with designing numerical schemes. As discussed in \cite{ar_Arora_1997,ar_Karni_1997} slowly moving shocks can generate large oscillations, and they propose the inclusion of additional diffusivity locally to the shock to eliminate these issues. We will discuss modifications to our (or indeed any) reconstruction that facilitate these local requirements for increased numerical diffusivity. We show that with piecewise constant reconstruction the central-upwind scheme \cite{ar_Kurganov_2001} is able to greatly suppress oscillations due to its increased numerical diffusivity. Thus, by using piecewise constant reconstruction around shocks, we achieve a greatly reduced wave train. In addition, the modification aids the simulation of dam break problems, in which fluid spreads over an initially dry bed, which removes the need for a tailwater (numerically imposing a non zero depth in the dry region). The modification is, by design, compatible with source terms, and therefore does not disrupt the well-balancing of the scheme.

The article is organised as follows. We begin with an overview of finite volume schemes in \cref{sec:FV}, followed by the reconstruction derived in \cite{ar_Skevington_N002_Recon_Theory} (\cref{sec:Recon}). In \cref{sec:SS} modifications to the reconstruction are made to ensure the resolution of slow shocks and dam break flows. The resulting algorithm is presented in \cref{sec:algorithm}. Other algorithms from the literature are overviewed in \cref{sec:others}, to which we will compare our new scheme. The new scheme is then used for a selection of numerical experiments (\cref{sec:TP}). We conclude in \cref{sec:END}.
\section{The numerical scheme}

\subsection{Semi-Discrete Finite-Volume Schemes}	\label{sec:FV}

Hyperbolic systems of conservation laws have the general form
\begin{equation}	\label{eqn:FV_conservation}
	\pdv{Q}{t} + \pdv{}{x}(F) = \Psi	\quad \textrm{for} \quad x_L \leq x \leq x_R,
\end{equation}
where $Q$ is a function from $(x,t)$ to $\mathbb{R}^M$, and the flux $F$ and source $\Psi$ are functions from $(Q,x,t)$ to $\mathbb{R}^M$. The system is spatially discretized over $J$ cells by introducing $J+1$ cell boundary points $x_L=x_{1/2}<x_{3/2}<\ldots<x_{J+1/2}=x_R$ of width $\Delta x_j=x_{j+1/2} - x_{j-1/2}$ with cell centres separated by a distance $\Delta x_{j+1/2} = (\Delta x_{j}+\Delta x_{j+1})/2$, and defining cell averages
\begin{align}\label{eqn:FV_QPdef} 
	Q_{j}(t) &= \frac{1}{\Delta x} \int_{x_{j-1/2}}^{x_{j+1/2}} Q \dd x,
	&\textrm{and}&&
	\Psi_{j}(t) &= \frac{1}{\Delta x} \int_{x_{j-1/2}}^{x_{j+1/2}} \Psi \dd x.
\end{align}
Performing the same averaging process over the entire system \cref{eqn:FV_conservation} yields
\begin{equation}
	\dv{Q_{j}}{t} + \frac{1}{\Delta x_j} \pbrk*{ F_{j+1/2} - F_{j-1/2} } = \Psi_{j}
\end{equation}
where $F_{j+1/2}(t) = F(Q(x_{j+1/2},t),x_{j+1/2},t)$. In order to construct a spatially discrete scheme we express $F_{j+1/2}$ and $\Psi_{j}$ approximately as functions of $Q_{j}$. We first perform a reconstruction of $Q$ within each cell, here we will use a piecewise linear reconstruction with the gradient in cell $j$ denoted by $[Q_x]_j$ (its computation discussed in \cref{sec:Recon}). Thus the limiting values at the cell boundaries are
\begin{align}	\label{eqn:FV_Qpm}
	Q_{j-1/2}^{+} &= Q_{j} - \frac{\Delta x_j}{2}[Q_x]_j,
	&
	Q_{j+1/2}^{-} &= Q_{j} + \frac{\Delta x_j}{2}[Q_x]_j,
\end{align}
\sloppy and the variation across a cell is then $\Delta Q_j \eqdef Q_{j+1/2}^{-} - Q_{j-1/2}^{+} = \Delta x_j [Q_x]_j$ whilst the variation between cells is $\Delta Q_{j+1/2} \eqdef Q_{j+1}-Q_j$, and similarly for other fields. From these we construct the approximation $F_{j+1/2}(t) \sim \hat{F}(Q_{j+1/2}^{-}(t),Q_{j+1/2}^{+}(t),x_{j+1/2},t)$, and we will employ the central-upwind scheme in our numerical tests \cite{ar_Kurganov_2001}. To evolve through time we use the second order SSP RK (Strong Stability-Preserving Runge Kutta) scheme from \cite{ar_Shu_1988,ar_Gottlieb_2001}. The boundary conditions are imposed using the approach detailed in \cite{ar_Skevington_N001_Boundary} (including the RKNR time-step modification) which employs forcing computed from first order extrapolation to ensure consistency with the solution in the bulk.

\subsection{Well balanced source and reconstruction}	\label{sec:Recon}

We begin with the expression of the geometric source term in \cref{eqn:INTRO_SW_Full_u}. The bed function we assume to be continuous, with discrete values at the cell interfaces $b_{j+1/2}$. The source is approximated as \cite{ar_Bermudez_1994}
\begin{equation}	\label{eqn:source_discretisation}
	-gh \dv{b}{x} \sim -g \frac{h_{j-1/2}^+ + h_{j+1/2}^-}{2} \frac{\Delta b_j}{\Delta x_j}
\end{equation}
which is known to be well balanced for a reconstruction with constant surface elevation \cite{ar_Bermudez_1994,ar_Kurganov_2001,ar_Skevington_N002_Recon_Theory}.

We use the reconstruction developed in \cite{ar_Skevington_N002_Recon_Theory}, modified in \cref{sec:SS} by including a discrete measure of the local smoothness $\Theta_j \in [0,1]$ which (at sufficiently high resolution) is $1$ where the solution is continuous and less than $1$ local to shocks. The reconstruction employs the $\minmod$ slope limiter \cite{ar_Tadmor_1988}
\begin{subequations}	\label{eqn:minmod_recon}
\begin{align}\label{eqn:minmod_recon_R}
	\sigma^\textrm{mm}(v_{j-1},v_{j},v_{j+1}) &= \minmod \pbrk*{
		\frac{2 \alpha_{j-1/2}^{+} (v_{j}-v_{j-1})}{\Delta x_j} , 
		\frac{2 \alpha_{j} (v_{j+1}-v_{j-1})}{\Delta x_j} , 
		\frac{2 \alpha_{j+1/2}^{-} (v_{j+1}-v_{j})}{\Delta x_j}
	},  
\end{align}
which is parametrised by $\alpha_j,\alpha_{j+1/2}^\pm \in [0,1]$, and
\begin{align}\label{eqn:minmod_recon_MM}
	\minmod[w_1 , w_2 \ldots w_n] \eqdef
	\begin{cases}
		\min[w_1 , w_2 \ldots w_n] 	& \textrm{if } \min[w_1 , w_2 \ldots w_n]>0,	\\
		\max[w_1 , w_2 \ldots w_n] 	& \textrm{if } \max[w_1 , w_2 \ldots w_n]<0,	\\
		0							& \textrm{otherwise}.
	\end{cases}
\end{align}
\end{subequations}
In our numerical tests we take $\alpha_{j+1/2}^{\pm} = 3/4$, $\alpha_{j} = 1/4$. To reconstruct the depth we first compute the two reconstructions, one in $h$ and one in $\eta$, as
\sidebysidesubequations{eqn:gradient_h_eta}
	{[h_x]_j^{h} \eqdef \Theta_j^h \sigma^\textrm{mm}(h_{j-1},h_{j},h_{j+1}),}
	{eqn:gradient_h}{}
	{[h_x]_j^{\eta} \eqdef \Theta_j^h \sigma^\textrm{mm}(\eta_{j-1},\eta_{j},\eta_{j+1}) - \frac{\Delta b_j}{\Delta x_j}.}
	{eqn:gradient_eta}
where $\Theta_j^h$ is similar to $\Theta_j$, see \cref{sec:algorithm}. From these gradients cell interface values may be computed and are denoted $h_{j+1/2}^{h\pm}$ and $h_{j+1/2}^{\eta\pm}$. We then compute the final gradient to be used as the convex combination
\begin{equation} \label{eqn:reaults_gradient_convex}
	[h_x]_j \eqdef (1-\gamma_j) [h_x]_j^{h} + \gamma_j [h_x]_j^{\eta}.
\end{equation}
from which the final reconstructed values are found and denoted $h_{j+1/2}^\pm$. The convex coefficient $\gamma_j$ has the expression
\begin{align}	\label{eqn:convex_coeff}
	\gamma_j(\xi_j) &\eqdef
	\begin{cases}
		0 & \textrm{if } \xi_j \leq 1,	\\
		G_j(\xi_j - 1)	& \textrm{if } 1 \leq \xi_j \leq \xi_j^{C},	\\
		1 & \textrm{if } \xi_j^{C} \leq \xi_j,
	\end{cases}
	&&\textrm{where}&
	\xi_j^{C} &\eqdef 1 + \frac{1}{G_j},
	&
	G_j &\eqdef 1-\max \pbrk*{ \alpha_{j-1/2}^{+} , \alpha_{j+1/2}^{-} },
\end{align}
and $\xi_j \eqdef {h_j^\downarrow}/{\Delta b_j^\uparrow}$ (when $\Delta b_j^\uparrow=0$ then $[h_x]_j$ is independent of $\xi_j$, so may set to any value) where
\begin{align}
	h_j^\downarrow &\eqdef \min \pbrk*{ h_j + \alpha_{j-1/2}^{+} \Delta h_{j-1/2} ,h_j,h_j + \alpha_{j+1/2}^{-} \Delta h_{j+1/2} },	
	\label{eqn:results_hmax}	\\
	\Delta b_j^\uparrow &\eqdef \max \pbrk*{ \abs*{\Delta b_j/2-\alpha_{j-1/2}^{+} \Delta b_{j-1/2}},\abs*{\Delta b_j/2},\abs*{\Delta b_j/2 - \alpha_{j}(b_{j+1}-b_{j-1})},\abs*{\Delta b_j/2 - \alpha_{j+1/2}^{-}  \Delta b_{j+1/2}} , B_j },
	\label{eqn:results_dbmax}\\
	B_j &\eqdef \ppar*{\frac{q_j^2}{\Fro_0^2 g}}^{1/3}.
	\label{eqn:B_fastfluid}
\end{align}
This ensures that the reconstruction is positivity preserving and bounded. In regions where the fluid is deep (relative to bed variation) and slow ($\Fro \ll \Fro_0$) then $\gamma_j = 1$ and the reconstruction is in $\eta$ alone, \ie $[h_x]_j = [h_x]_j^\eta$. In regions were the fluid is shallow (relative to bed variation) or fast ($\Fro \gg \Fro_0$) such as in thin film regions, then $\gamma_j = 0$ and the reconstruction is in $h$ alone, \ie $[h_x]_j = [h_x]_j^h$. For our simulations we employ $\Fro_0 = 10 \approx (\xi_j^C)^{3/2}$. The inclusion of the suppressor $\Theta_j$ in the depth reconstruction \cref{eqn:gradient_h_eta} was not done in \cite{ar_Skevington_N002_Recon_Theory} and is a new addition here to be discussed in \cref{sec:SS}.

The gradient of flux is computed using the reconstruction
\begin{equation}\label{eqn:gradient_q}
	[q_x] \eqdef \Theta_j^q \sigma^{\textrm{mm}} (q_{j-1} , q_j , q_{j+1})
\end{equation}
which yields a bounded velocity by the arguments in \cite[\S\S4,8]{ar_Skevington_N002_Recon_Theory}. Here $\Theta_j^q$ is similar to $\Theta_j$, see \cref{sec:algorithm}.

\subsection{Non-oscillatory slow shocks}	\label{sec:SS}

Many solutions to the shallow water equations, and hyperbolic systems in general, contain discontinuous shocks. When a compressive shock (\ie one where the characteristics converge on the shock, as opposed to a contact shock where they are locally parallel) exists in a numerical scheme it may generate wave pulses at regular intervals. For slowly moving shocks in the shallow water equations this occurs as a result of the spike in volume flux \cite{ar_Montilla_2017}. This is a large value of $q$ that exists within a shock because of the shape of the Hugoniot locus in $(h,q)$ space. That is to say, the large value of $q$ is evidence that the numerical scheme is approximately resolving the shock conditions. As the shock moves, this spike moves with it, and as it passes from one cell to the next a wave pulse is generated, as can be seen in many numerical experiments \cite{ar_Lin_1995,ar_Jin_1996,ar_Karni_1997,ar_Arora_1997,ar_Striba_2003,ar_Johnsen_2013,ar_Montilla_2017}. For some first order schemes it is possible to eliminate the spike by using specialist techniques \cite{ar_Montilla_2017}. For higher order schemes no such strategy is known, instead additional numerical diffusion may be added to widen the spike, which can dramatically reduce the size of the oscillations generated as it moves \cite{ar_Lin_1995,ar_Karni_1997,ar_Arora_1997}. Our numerical experiments have found that, for the central-upwind scheme \cite{ar_Kurganov_2001}, the numerical diffusion that results from piecewise constant reconstruction is sufficient to suppress the oscillations (see \cref{tp:TP_slowshock}). One approach that may be used employs the measure from \cite[(2.9)]{ar_Skevington_N002_Recon_Theory}, specifically taking $\Theta_j^q = \kappa_j$  where
\begin{align}	\label{eqn:basic_suppressor}
	\kappa_j \eqdef \min \pbrk*{ 1 , \frac{K_{j-1/2}^+ h_j}{h_{j-1}} , \frac{K_{j+1/2}^- h_j}{h_{j+1}} }
\end{align}
($1/0 = +\infty$ so that the $\min$ does not select this value, $0/0=0$ to that, if this cell is dry, the gradient is $0$). For $K_{j+1/2}^\pm = 1 + \order{1}$ as $\Delta x \rightarrow 0$ this choice may have the desired properties local to shocks (see \cref{tp:TP_slowshock}). However, it also detects regions where the solution transitions from lake at rest to dry, as will any detector based on the magnitude of derivatives relative to a local depth measure. Thus, if wish to obtain an accurate solution local to the edge of a lake then we require $K_{j+1/2}^\pm \gg 1$, which means that it will not detect most shocks.

To construct a shock detector that does not detect regions where the solution dries, we must be precise about what we mean by a shock. To do this we will use fundamental properties of the general conservation law \cref{eqn:FV_conservation}. In particular, that the matrix $\pdv*{F}{Q}$ has real eigenvalues $\lambda^{(m)}$ corresponding to left eigenvectors $l^{(m)}$, with $m$ indexing over the characteristic fields. Left multiplying \cref{eqn:FV_conservation} by $l^{(m)}$ results in the characteristic equation for the $m^\nth$ filed
\begin{align}
	l^{(m)} \ppar*{ \pdv{Q}{t} + \lambda^{(m)} \pdv{Q}{x} } &= l^{(m)} \Psi'
	&&\textrm{where}&
	\Psi' \eqdef \Psi - \pdv{F}{x}.
\end{align}
By the Lax entropy condition \cite{ar_Lax_1957,bk_Lax_SP1} characteristics converge on a discontinuity. That is for a discontinuity at location $x_s(t)$, defining $\lambda_-^{(m)}$ to be a characteristic speed to the left of the shock and $\lambda_+^{(m)}$ to the right, there must be some field $m$ so that $\lambda_-^{(m)} \geq \dv*{x_s}{t} \geq \lambda_+^{(m)}$. We distinguish between contact discontinuities with $\lambda_-^{(m)} = \lambda_+^{(m)}$ and shocks for which $\lambda_-^{(m)} > \lambda_+^{(m)}$. Numerically, contact discontinuities are diffused over some region and are not sharpened by the action of fluxes, thus we do not wish to increase the diffusivity local to these. Thus we look for singular gradients in the characteristic speeds $\lambda^{(m)}$, that is $\lambda^{(m)}(x-\epsilon,t) - \lambda^{(m)}(x+\epsilon,t)$ tends to some finite positive quantity as $\epsilon \rightarrow 0$. We define the quantity
\begin{equation}	\label{eqn:SS_detector_speed}
	\Delta \lambda_j^{(m)} \eqdef \frac{x_\textrm{ref}^{p_1}}{\lambda_\textrm{ref}} \max \pbrk*{
		\frac{\lambda_{j-1}^{(m)} - \lambda_{j}^{(m)}}{(\Delta x_{j-1/2})^{p_1}} , 
		\frac{\lambda_{j}^{(m)} - \lambda_{j+1}^{(m)}}{(\Delta x_{j+1/2})^{p_1}} , 
		0 }
\end{equation}
where $\lambda_\textrm{ref}$ is the speed-scale and $x_\textrm{ref}$ the length-scale so that $\Delta \lambda_j^{(m)}$ is dimensionless (for our simulations of dimensionless systems we take $\lambda_\textrm{ref}=1$ and $x_\textrm{ref}=x_R-x_L$) and $\lambda_j^{(m)} \eqdef \lambda^{(m)}(Q_j,x_j,t)$. This expression is $\order{\Delta x^{1-p_1}}$ in smooth regions and $\order{\Delta x^{-p_1}}$ local to shocks, so taking $0<p_1<1$ means that \cref{eqn:SS_detector_speed} can be used as a detector as $\Delta x \rightarrow 0$.

Whilst \cref{eqn:SS_detector_speed} is a reliable detector of shocks, it suffers from some of the same problems as \cref{eqn:basic_suppressor}. In particular, with source terms it is possible to construct situations where characteristics are converging rapidly, and yet the dominant balance in the system is between the flux gradient and the source. This happens in the shallow water system around the edge of a quiescent lake, \cref{fig:PWB_draining_flow}, especially when the draining layer is orders of magnitude shallower than the lake (there is, formally, a shock between the draining layer and the lake, but in the simulation there is insufficient resolution to capture it). To detect this situation we define
\begin{align}	\label{eqn:SS_detector_flux}	\begin{split}
	\Delta F_j^{(m)} \eqdef x_\textrm{ref}^{p_1-1}\max   & \left[
		\abs*{ \frac**{\frac{\lambda_{j-1}^{(m)} l_{j-1}^{(m)} \Delta Q_{j-1/2} 	}{ (\Delta x_{j-1/2})^{p_1} }}		{l_{j-1}^{(m)} \Psi'_{j-1}} } ,
		\abs*{ \frac**{\frac{\lambda_{j}^{(m)} l_{j}^{(m)} \Delta Q_{j-1/2} 		}{ (\Delta x_{j-1/2})^{p_1} }}		{l_{j}^{(m)} \Psi'_{j}}     } ,
		\right. \\& \left.
		\abs*{ \frac**{\frac{\lambda_{j}^{(m)} l_{j}^{(m)} \Delta Q_{j+1/2}		}{ (\Delta x_{j+1/2})^{p_1} }}		{l_{j}^{(m)} \Psi'_{j}}     } ,
		\abs*{ \frac**{\frac{\lambda_{j+1}^{(m)} l_{j+1}^{(m)} \Delta Q_{j+1/2}	}{ (\Delta x_{j+1/2})^{p_1} }}		{l_{j+1}^{(m)} \Psi'_{j+1}} } 
		\right].
\end{split}\end{align}
which is the dimensionless ratio of the strength of flux gradient (up to a power of $\Delta x$) to source terms in each characteristic field. This expression is $\order{\Delta x^{1-p_1}}$ in smooth regions, and $\order{\Delta x^{-p_1}}$ local to shocks where the flux gradient dominates the source.

We say that there is a shock in the $m^\nth$ field when both $\Delta \lambda_j^{(m)} \gg 1$ and $\Delta F_j^{(m)} \gg 1$. We define a shock detector for the $m^\textrm{th}$ field as the dimensionless quantity
\begin{equation}	\label{eqn:SS_detector_field}
	\Theta_j^{(m)} = 1 - \ppar*{ \frac{1}{\ppar*{\Delta \lambda_j^{(m)}}^{p_2}} + 1 }^{-p_3} \ppar*{ \frac{1}{\ppar*{\Delta F_j^{(m)}}^{p_2}} + 1 }^{-p_3}.
\end{equation}
In smooth regions $\Delta \lambda_j^{(m)} \ll 1$ and $\Delta F_j^{(m)} \ll 1$, thus
\begin{subequations}
\begin{align}
	\Theta_j^{(m)} &\sim 1 - \ppar*{\Delta \lambda_j^{(m)}}^{p_2 p_3} \ppar*{\Delta F_j^{(m)}}^{p_2 p_3} = 1 - \order{ \Delta x^{2 p_2 p_3 (1-p_1)}},
\intertext{whilst around contact discontinuities $\Delta \lambda_j^{(m)} \ll 1$ and $\Delta F_j^{(m)} \gg 1$, thus}
	\Theta_j^{(m)} &\sim 1 - \ppar*{\Delta \lambda_j^{(m)}}^{p_2 p_3} = 1 - \order{ \Delta x^{p_2 p_3 (1-p_1)}},
\intertext{and if characteristics are converging on a location where the dominant balance is between flux and source then $\Delta \lambda_j^{(m)} \gg 1$ and $\Delta F_j^{(m)} \ll 1$, thus}
	\Theta_j^{(m)} &\sim 1 - \ppar*{\Delta F_j^{(m)}}^{p_2 p_3} = 1 - \order{ \Delta x^{p_2 p_3 (1-p_1)}},
\intertext{finally in regions where there are shocks $\Delta \lambda_j^{(m)} \gg 1$ and $\Delta F_j^{(m)} \gg 1$ thus}
	\Theta_j^{(m)} &\sim p_3 \ppar*{ \ppar*{\Delta \lambda_j^{(m)}}^{-p_2} + \ppar*{\Delta F_j^{(m)}}^{-p_2} } = \order{ \Delta x^{p_2 p_1}}.
\end{align}
\end{subequations}
To obtain a second order scheme we require $p_2 p_3 (1-p_1) \geq 1$, and to suppress the gradients local to shocks $p_2 p_1 > 0$, thus $0 < p_1 p_2 \leq p_2 - (1/p_3)$. The values $p_1 = 1/2$, $p_2 = 2$, $p_3=1$ are appropriate for our scheme, so that $\Theta_j^{(m)} = 1 - \order{ \Delta x^{2} }$ in smooth regions. For higher order reconstruction the value of $p_3$ may be increased so that $1-\Theta_j^{(m)} = \order{ \Delta x^{2 p_3} }$ is smaller.

To ensure bounded velocities we additionally employ a generalisation of the measure from \cite{ar_Skevington_N002_Recon_Theory}, \cref{eqn:basic_suppressor}, to all positive fields $m \in \mathscr{P}$
\begin{align}	\label{eqn:SS_detector_dry}
	\hat{\Theta}_j^{(m)} &= \min \pbrk*{ 1 ,\ppar*{ K_{j-1/2}^+ \frac{Q_j^{(m)}}{Q_{j-1}^{(m)}} }^{p_4} , \ppar*{ K_{j+1/2}^- \frac{Q_j^{(m)}}{Q_{j+1}^{(m)}} }^{p_4} }.
\end{align}
where $K_{j+1/2}^\pm \gg 1$ to ensure a smooth transition to dry cells, and the notation $Q^{(m)}$ indicates the $m^\nth$ component of $Q$. The large value of $K_{j+1/2}^\pm$ means that the velocity bound is also large, but for $p_4 \geq 1$ this large bound only exists for ${h_{j \pm 1}^\mp}/{h_j} \approx K_{j \pm 1/2} ^{\mp}$, else the bound can be substantially reduced. For our scheme $K_{j \pm 1/2} ^{\pm} = 100$ and $p_4=2$ are appropriate.

Finally we define
\begin{equation}	\label{eqn:SS_detector_compile}
	\Theta_j = \min \pbrk*{ \min_m \Theta_j^{(m)} , \min_{m \in \mathscr{P}} \hat{\Theta}_j^{(m)} }.
\end{equation}
This identifies shocks in any of the characteristic fields in such a way that the transition to dry cells is accurately resolved. Our reconstruction is no longer formally self monotone \cite{ar_Skevington_N002_Recon_Theory}, but this issue is only present when $1-\Theta_j = \order{1}$, \ie local to strong shocks or when the depth field is orders of magnitude shallower than its neighbours, in which case suppressing oscillations and ensuring bounded velocities (respectively) are more important for stability.

While this method was originally developed and presented here with consideration of slow shocks, it does not depend on the speed of the shock, or even the existence of a shock at all, only the convergence of characteristics and the domination of source terms by the flux gradient. In many problems characteristics converge at a wetting front, such as occurs in a dam break flow. Our method greatly aids the simulation of such phenomena, as will be demonstrated later.

\subsection{Implementation} \label{sec:algorithm}

We close this section by presenting our reconstruction as an algorithm to be implemented. We assume that the discrete fields and lattice points are known at the time $t$, and proceed as follows.
\begin{enumerate}
	\item 	Compute the source using the approximation \cref{eqn:source_discretisation}
	\item	Compute eigenvalues and eigenvectors from the cell centred values as $\lambda_j^{(m)} \eqdef \lambda^{(m)}(Q_j , x_j , t)$, $l_j^{(m)} \eqdef l^{(m)}(Q_j , x_j , t)$.
	\item 	Compute $\Theta_j^h = \Theta_j^q = \Theta_j$ using \cref{eqn:SS_detector_speed,eqn:SS_detector_flux,eqn:SS_detector_field,eqn:SS_detector_dry,eqn:SS_detector_compile}
	\item 	Perform the gradient reconstructions \cref{eqn:gradient_h_eta,eqn:gradient_q}
	\item	Calculate the convex coefficient \cref{eqn:convex_coeff} using \cref{eqn:results_hmax,eqn:results_dbmax,eqn:B_fastfluid}
	\item 	Assemble the final depth gradient \cref{eqn:reaults_gradient_convex}
\end{enumerate}
This yields the gradients $[h_x]_j$ and $[q_x]_j$ to be used to construct the values at the cell interfaces, and thereby the numerical flux between cells. The scheme presented here we will call \emph{SkT}, in the following section we briefly introduce the other schemes that we will compare it to.
\section{Other schemes} \label{sec:others}

In our numerical tests we will compare our scheme to a number of others developed in the literature. While only the scheme from \cite{ar_Skevington_N002_Recon_Theory} and its modification here were designed using convex combination \cref{eqn:gradient_h_eta,eqn:reaults_gradient_convex}, several other schemes in the literature may be reinterpreted in this framework, which is how they will be presented below. Thus all that remains to present is expressions for $\gamma_j$, $\Theta_j^h$ and $\Theta_j^q$. For all numerical tests we will compare well-balanced schemes based off of the $\minmod$ reconstruction.

\paragraph{Kurganov and Levy, 2002 \cite{ar_Kurganov_2002}} Here the reconstruction is in $\eta$ so long as the depth in the current cell and its neighbours is above some threshold $H$, otherwise the reconstruction is in $h$, that is
\begin{equation}
	\gamma_j = 
	\begin{cases}
		0	&	\text{if } \min[h_{j-1} , h_{j} , h_{j+1}] < H,	\\
		1	&	\text{otherwise}.
	\end{cases}
\end{equation}
We will employ $H = 0.1$. There is no gradient suppression to piecewise constant, $\Theta_j^{h} = \Theta_j^{q} = 1$. We term this scheme \emph{Ku01}.

\paragraph{Kurganov and Petrova, 2007 \cite{ar_Kurganov_2007}} Here the reconstruction is in $\eta$, but then altered to set the depth to precisely zero whenever it would otherwise be negative, and adjusting the gradient in the cell and the reconstruction at the other end of the cell accordingly. This may be expressed the the framework of convex combination as
\begin{equation}
	\gamma_j = 
	\begin{cases}
		\frac**{h_{j-1/2}^h}{\ppar*{h_{j-1/2}^h - h_{j-1/2}^\eta}}		&	\text{if } h_{j-1/2}^{\eta+} < 0,	\\
		\frac**{h_{j+1/2}^h}{\ppar*{h_{j+1/2}^h - h_{j+1/2}^\eta}}		&	\text{if } h_{j+1/2}^{\eta-} < 0,	\\
		1	&	\text{otherwise}.
	\end{cases}
\end{equation}
There is no gradient suppression, $\Theta_j^{h} = \Theta_j^{q} = 1$, though after reconstruction the fluxes are corrected to ensure finite values of the reconstructed velocities
\begin{equation}
	q_{j+1/2}^{\pm \text{corr}} 
	= q_{j+1/2}^\pm \cdot \sqrt{\frac{2}{1+\max\pbrk*{1,(\epsilon/h_{j+1/2}^\pm)^4}}}
	= q_{j+1/2}^\pm \cdot \sqrt{\frac{2 (h_{j+1/2}^\pm)^4}{(h_{j+1/2}^\pm)^4+\max\pbrk*{(h_{j+1/2}^\pm)^4,\epsilon^4}}}
\end{equation}
where $\epsilon = \Delta x$ is the desingularization parameter. We term this scheme \emph{Ku07}.

\paragraph{Chertock \etal, 2015 \cite{ar_Chertock_2015}} Similar to above the reconstruction is in $\eta$ unless it would generate negative depths, in which case the gradient is set to zero. This can be expressed as 
\begin{equation}
	\gamma_j = 
	\begin{cases}
		\frac**{[h_x]_j^h}{\ppar*{[h_x]_j^h - [h_x]_j^\eta}}			&	\text{if } h_{j-1/2}^{\eta+} < 0 \text{ or } h_{j+1/2}^{\eta-} < 0,	\\
		1	&	\text{otherwise},
	\end{cases}
\end{equation}
and $\Theta_j^h = 1$. Instead of reconstructing in flux $q$, the velocity $u$ is reconstructed and then the flux reconstruction is computed as the product of the velocity and depth. Note that this means the scheme is not formally a finite volume scheme as the reconstructed $q$ does not satisfy \cref{eqn:FV_QPdef}. Prior to reconstruction the velocity is desingularized using
\begin{equation}  \label{eqn:desing_Chertock}
	u_j^\text{corr} 
	= \frac{q_j}{h_j} \cdot \frac{2}{1+\max\pbrk*{1,(\epsilon/h_j)^2}}
	= q_j \cdot \frac{2 h_j}{h_j^2+\max\pbrk*{h_j^2,\epsilon^2}}
\end{equation}
where $\epsilon = 10^{-8}$. The velocity is then reconstructed using the minmod slope limiter \cref{eqn:minmod_recon_R}. We term this scheme \emph{Ch15}.

\paragraph{Skevington, YEAR \cite{ar_Skevington_N002_Recon_Theory}} We also present results for the scheme without the modifications in \cref{sec:SS}, that is we use \cref{eqn:convex_coeff}, $\Theta_j^h = 1$ and $\Theta_j^q = \kappa_j$ \cref{eqn:basic_suppressor}. In our numerical tests we use $K_{j-1/2}^+ = K_{j-1/2}^+ = 1 + 10 \Delta x_j / (x_R-x_L)$. We term this scheme \emph{SkK}.

The final scheme tested is not able to be written in terms of a convex combination.

\paragraph{Bollermann \etal, 2013} This scheme uses a reconstruction in $\eta$ everywhere except in cells that are deemed to be at the edge of a lake, or are too shallow to be reconstructed reliably, in which case multiple pieces may be used in the cell to better approximate the surface elevation. Similarly to Ch15, reconstruction is not performed in $q$ but rather $u$, being desingularized as
\begin{equation} \label{eqn:desing_Boller}
	u_j^\text{corr} 
	=  
	\begin{cases}
		\frac*{q_j}{h_j}		& \text{if } h_j \geq \epsilon,	\\
		0					& \text{otherwise},
	\end{cases}
\end{equation}
with $\epsilon = 10^{-9}$. We term this scheme \emph{Bo13}.

\section{Test Problems} \label{sec:TP}

We simulate a selection of test problems for the shallow water system \eqref{eqn:INTRO_SW_Full} made dimensionless so that $g=1$. The test problems will be simulated on uniform grids at resolutions of
\begin{equation}
	J \in \setpred*{ \round*{ 10^{a/4} } }{ a \in \mathbb{N} \textrm{ and } 8 \leq a \leq 16 }  = \{ 100,178,316,562,1000,1778,3162,5623,10000 \}.
\end{equation}
For our convergence analysis we will use $\ell_1$ error, which is computed by taking the absolute difference between simulated and exact values for every variable and location stated, and then averaging over all.

The first three problems focus on our reconstruction of the transition between wet and dry states. In \cref{tp:TP_lakerest} we examine the accuracy with which the scheme resolves a steady state  by considering an initially quiescent fluid at rest in a basin. This is then modified by imposing a thin film of fluid where there was a dry slope previously, so that in \cref{tp:TP_lakedrain} we examine the accuracy with which the scheme drains a slope of fluid. For \cref{tp:TP_lakeoscillate} we consider an oscillating fluid in a parabolic basin so that we can illustrate how our scheme deals with a moving contact point. For our final two problems we move away from bed topography to examine gradient suppression, firstly for slow shocks in \cref{tp:TP_slowshock} where we show that oscillations are greatly reduced, and secondly for a dam break over a flat bed for \cref{tp:TP_dambreak} demonstrating the ability of our scheme to resolve a dam break without requiring a tailwater.

\begin{figure}[tp!]
	\begin{center}
		\includegraphics{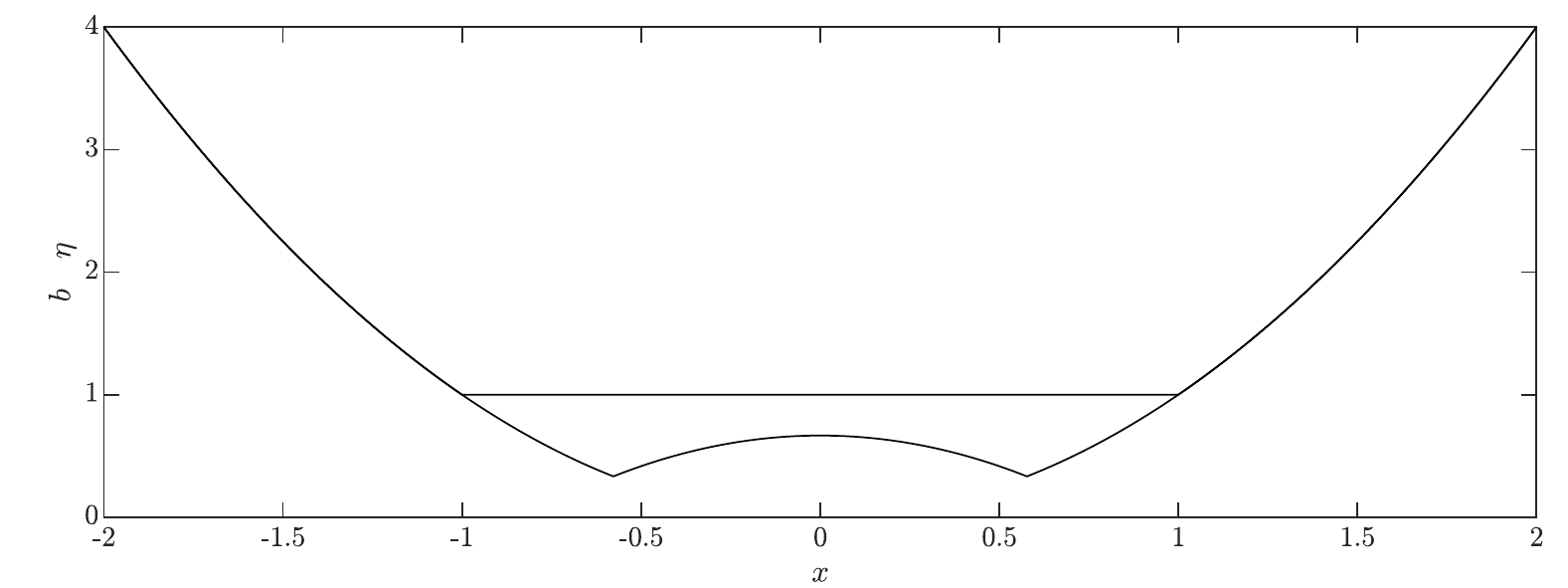}
	\end{center}
	\caption{The surface elevation $\eta$ and bed profile $b$ as functions of $x$ for \cref{tp:TP_lakerest}, taken from a simulation using SkT at resolution $J = 10^4$ at time $t=100$.}
	\label{fig:TP_lakerest}
\end{figure}
\begin{figure}[tp!]
	\begin{center}
		\includegraphics{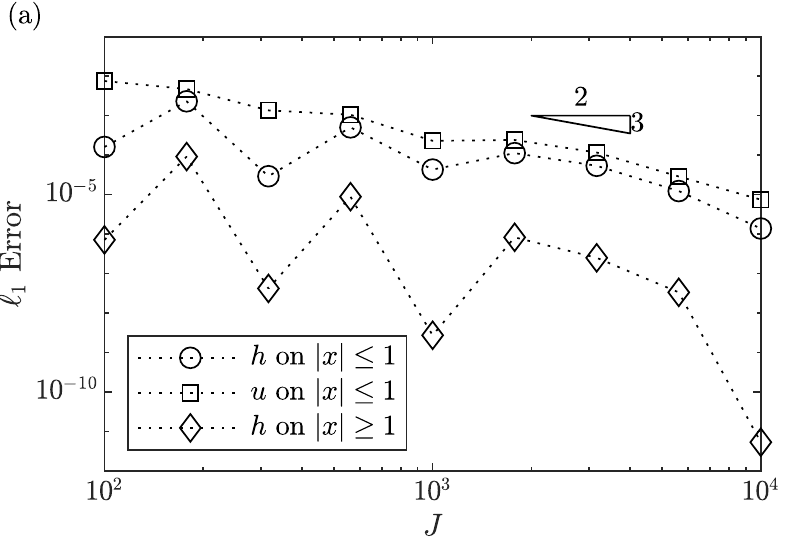}
		\includegraphics{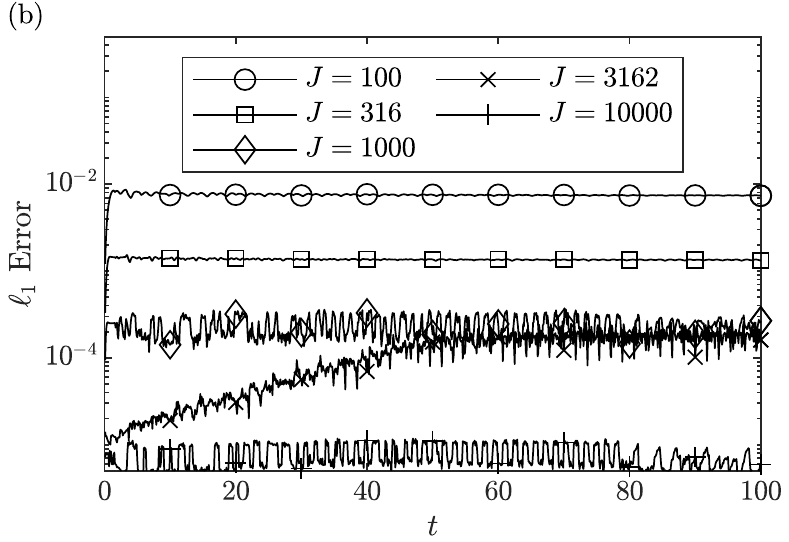}
	\end{center}
	\caption{The $\ell_1$ error between the simulations using SkT and the exact solution as a function of the resolution $J$ for \cref{tp:TP_lakerest}. In (a) we plot the error averaged over time $0 \leq t \leq 100$, as a function of resolution, separating the error between the wet ($\abs{x} \leq 1$) and dry ($\abs{x} \geq 1$) regions. In (b) we plot the error in $u$ in the wet region as a function of time. }
	\label{fig:TP_lakerest_convergence}
\end{figure}
\begin{figure}[tp!]
	\begin{center}
		\includegraphics{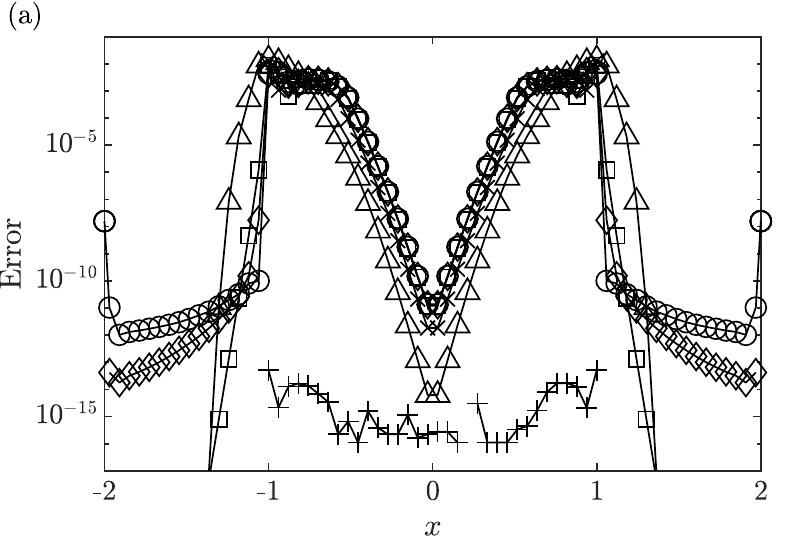}
		\includegraphics{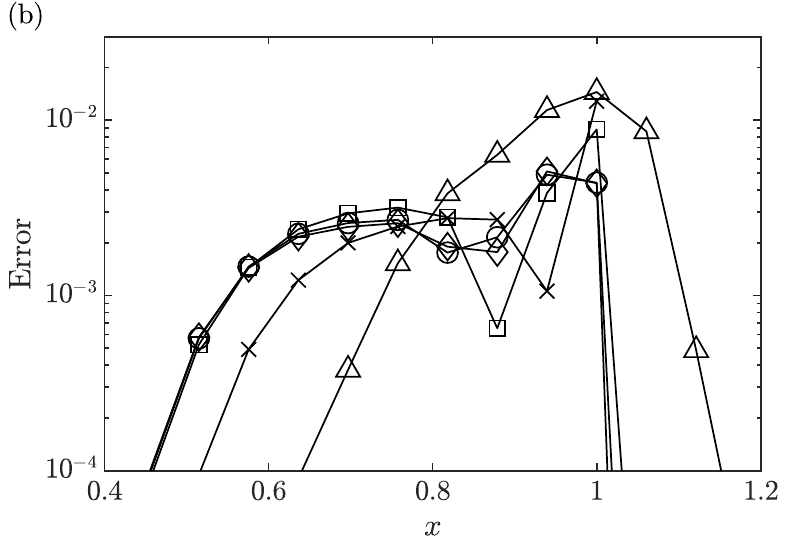}
	\end{center}
	\caption{The error in the depth field at $t=0.5$ for \cref{tp:TP_lakerest} as a function of $x$. The simulations were performed at a resolution of $J=66$. Different markers indicate different reconstructions: $\Diamond$ is Ku02, $\times$ is Ku07, $+$ is Bo13, $\Delta$ is Ch15, $\bigcirc$ is SkK, $\square$ is SkT. Plot (b) is a zoom in of plot (a).}
	\label{fig:TP_lakerest_comparison}
\end{figure}

\begin{testproblem}[Quiescent fluid in a parabolic basin]	\label{tp:TP_lakerest}
We begin by simulating a steady state, specifically a parabolic basin with a bump $b = \abs{x^2-1/3}+1/3$ in which there is fluid with surface $\eta=\max(1,b)$ and velocity $u=0$, see \cref{fig:TP_lakerest}. We take this as an initial condition for our simulation, which we perform on the interval $-2 \leq x \leq 2$ with boundary conditions $u=0$ for $x \in \{-2,2\}$. \Cref{fig:TP_lakerest_convergence} indicates convergence around $\order{J^{-3/2}}$ with some variation depending on the position of the contact line within the cell. In addition we plot the error as a function of time, showing that the error is independent of the duration of the simulation.

Comparing the different reconstructions (\cref{fig:TP_lakerest_comparison}), Bo13 resolves at machine precision, and so if a pure steady state is required then this is the one to choose. Among the others, Ku07 has a reasonably small error at the edge of the lake, and no fluid on the slope, and SkT has a similar error across the domain, and  rapidly decaying error on the slope. The Ku02 simulation failed at time $t=0.78$ (2 \sf) due to the problem discussed in \cite[fig. 3.1]{ar_Skevington_N002_Recon_Theory}, to remedy this would require ensuring bounded velocities.
\end{testproblem}

\begin{figure}[tp!]
	\begin{center}
		\includegraphics{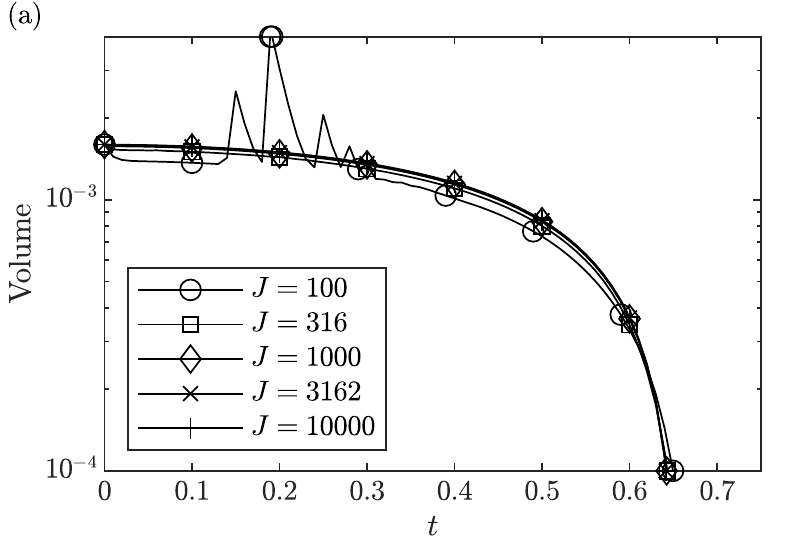}
		\includegraphics{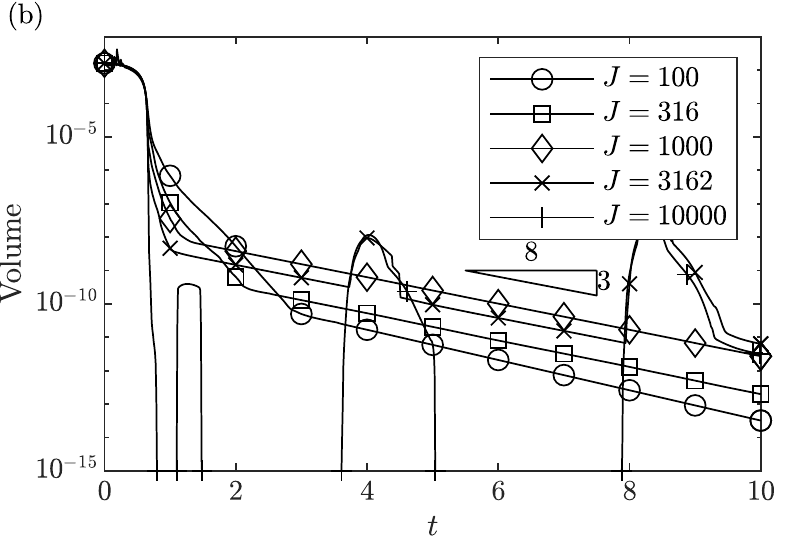}
	\end{center}
	\caption{The volume of fluid in $\abs{x} \geq 1.2$ as a function of $t$ for \cref{tp:TP_lakedrain} simulated using SkT. In (a) we plot early times where in the exact solution we have some fluid, and (b) later times.}
	\label{fig:TP_lakedrain_convergence}
\end{figure}
\begin{figure}[tp!]
	\begin{center}
		\includegraphics{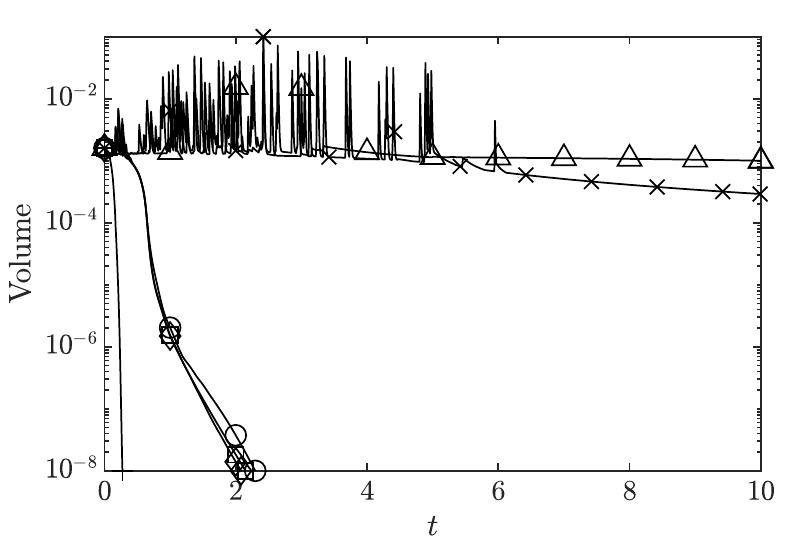}
	\end{center}
	\caption{The volume of fluid in $\abs{x} \geq 1.2$ as a function of $t$ for \cref{tp:TP_lakedrain}. The simulates were performed at a resolution of $J=66$. Different markers indicate different reconstructions: $\bigcirc$ is SkK, $\square$ is SkT, $\Diamond$ is Ku02, $\times$ is Ku07, $+$ is Bo13,  $\Delta$ is Ch15.}
	\label{fig:TP_lakedrain_comparison}
\end{figure}

\begin{testproblem}[Draining into a parabolic basin]	\label{tp:TP_lakedrain}
To investigate the numerical draining of fluid down a slope we modify \cref{tp:TP_lakerest} by setting $\eta=\max(1,b+10^{-3})$ and $u=0$ at time $t=0$. The initial condition appears very similar to \ref{fig:TP_lakerest}, with a lake of fluid in $-1 < x < 1$, and a thin film of fluid on the slopes $x < -1$ and $1 < x$. An asymptotic solution the this problem can be constructed. The depth on the slopes is very small, thus we take $h=\epsilon \tilde{h}$, the momentum equation is then, to leading order
\begin{align}
	\pdv{u}{t} + u \pdv{u}{x} &\sim -\pdv{b}{x},
	&
	\therefore\qquad
	\dv[2]{x}{t} &\sim -2x \quad \textrm{on} \quad \dv{x}{t} = u
\end{align}
which are the trajectories of fluid parcels. From here we arrive at the solution $x=x_0 \cos (\sqrt{2} t)$ (because $u=0$ at $t=0$), thus the fluid that started at $x=2$ (the edge of the domain) reaches $x=1$ (the lake) at time $t=\frac*{\pi}{(2 \sqrt{2})} = 0.74$ (2 \sf). Beyond this time, aside from the effect of small ripples on the lake initiated by the draining fluid, the slopes should be dry (\ie $h=0$ on $\abs{x} \gtrsim 1$). To avoid problems associated with the ripples we plot the total volume of fluid in the region $\abs{x} \geq 1.2$ in \cref{fig:TP_lakedrain_convergence}; this region is dry by $t = 0.66$ (2 \sf) in the asymptotic solution, a fact with which the simulations are consistent. Beyond this time we observe exponential convergence in time to the dry bed steady state with the volume on the slope at order $10\hat{\;\;}(-3t/8)$, the exponent determined approximately from the simulation results.

Comparing the different reconstructions (\cref{fig:TP_lakedrain_comparison}), only SkK, SkT, and Ku02 accurately capture the draining of the slope. For the others, Ku07 and Ch15 do not drain, while Bo13 drains faster than the exact solution. The latter is because, on draining slopes, Bo13 sets the upslope depth to zero and the downslope depth to be significantly greater than the cell average depth, and as a consequence the diffusive numerical error across each cell interface drives the fluid down the slope and into the lake. 
\end{testproblem}

\begin{figure}[tp!]
	\begin{center}
		\includegraphics{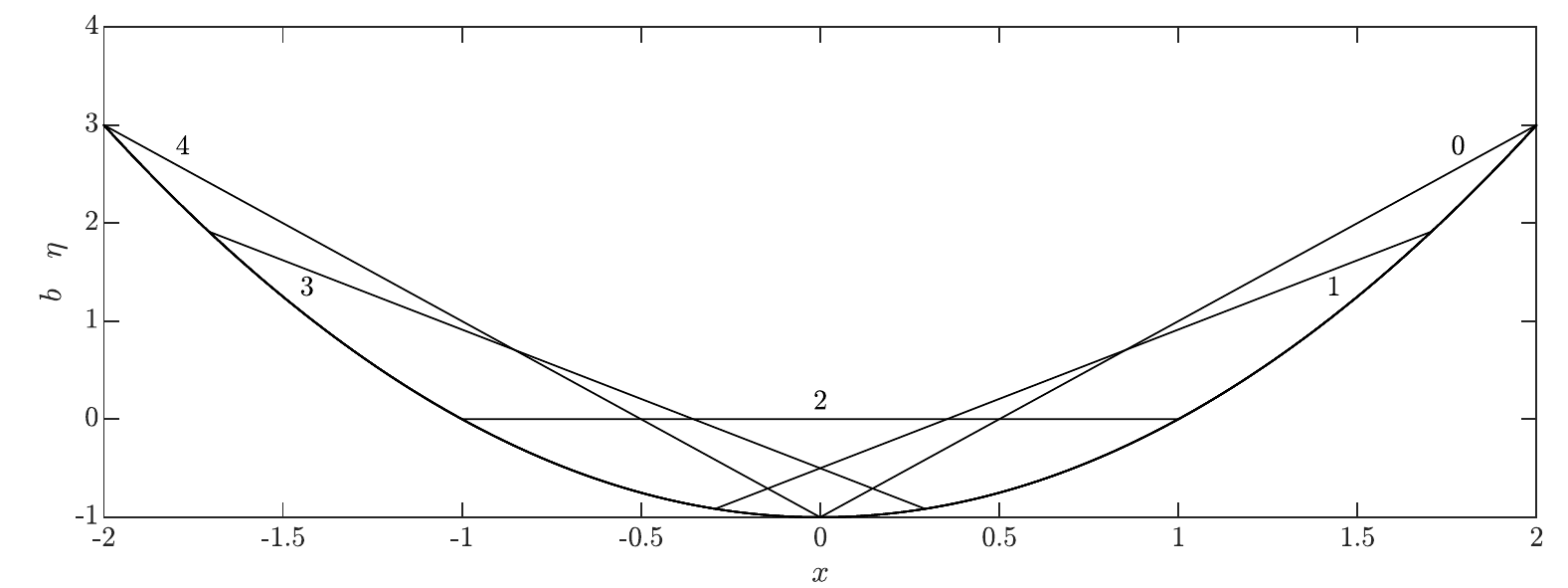}
	\end{center}
	\caption{The depth field as a function of $x$ for \cref{tp:TP_lakeoscillate} at times ${8t}/{(\sqrt{2}\pi)} = 0,1,2,3,4$ (period of oscillation is $\sqrt{2}\pi$), these values annotate the plot. The data is taken from a simulation using SkT at resolution $J = 10^4$.}
	\label{fig:TP_lakeoscillate}
\end{figure}
\begin{figure}[tp!]
	\begin{center}
		\includegraphics{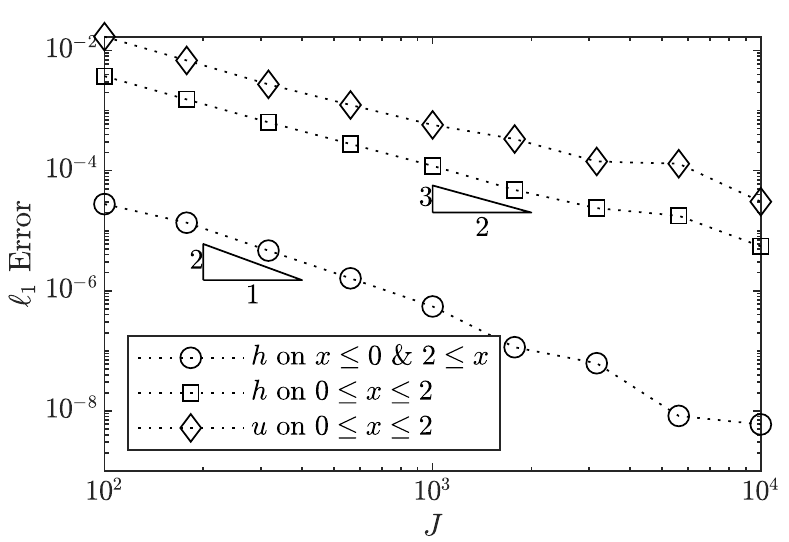}
	\end{center}
	\caption{The $\ell_1$ error between the simulations using SkT and the exact solution as a function of the resolution $J$. The error is evaluated at time $t=\sqrt{2}\pi$, and separately for the wet and dry regions.}
	\label{fig:TP_lakeoscillate_convergence}
\end{figure}
\begin{figure}[tp!]
	\begin{center}
		\includegraphics{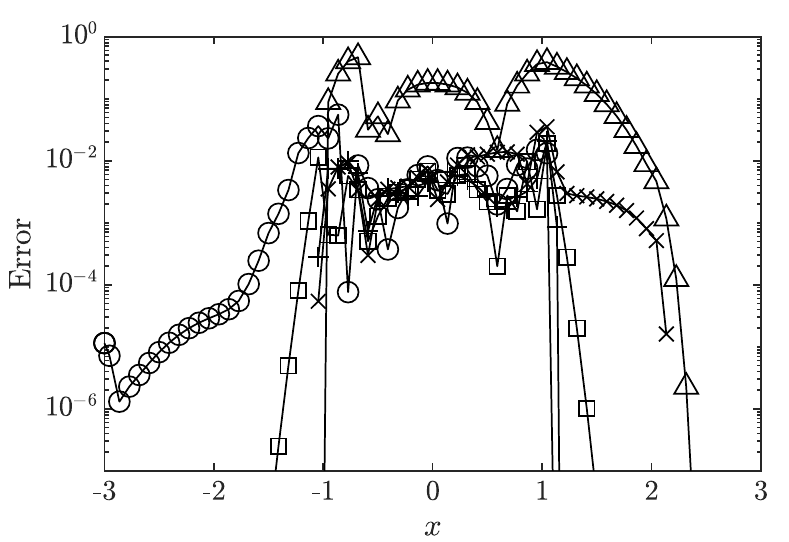}
	\end{center}
	\caption{The error in the depth field at $t=\sqrt{2}\pi/4$ for \cref{tp:TP_lakeoscillate} as a function of $x$. The simulations were performed at a resolution of $J=66$. Different markers indicate different reconstructions: $\Diamond$ is Ku02, $\times$ is Ku07, $+$ is Bo13, $\Delta$ is Ch15, $\bigcirc$ is SkK, $\square$ is SkT.}
	\label{fig:TP_lakeoscillate_comparison}
\end{figure}

\begin{testproblem}[Oscillations in a parabolic basin]	\label{tp:TP_lakeoscillate}
In this final test problem using bed topography we discuss the planar oscillation in a parabola that was proposed by Thacker \cite{ar_Thacker_1981}, see \cite{ar_Delestre_2016}. Specifically, we consider the dynamics on the topography $b(x) = x^2-1$ with solution (see \cref{fig:TP_lakeoscillate})
\begin{align}
	h &=
	\begin{cases}
		1 - (x-\cos(\sqrt{2}t))^2 	& \text{for } -1 \leq x-\cos(\sqrt{2}t) \leq 1,\\
		0 							& \text{otherwise},
	\end{cases}
	&
	u &=
	\begin{cases}
		-\sqrt{2} \sin(\sqrt{2}t) 	& \textrm{for } -1 \leq x-\cos(\sqrt{2}t) \leq 1,\\
		\text{undefined} 			& \text{otherwise}.
	\end{cases}
\end{align}
The depth field is plotted in \cref{fig:TP_lakeoscillate}. Setting $B_j = 0$ everywhere causes the simulation to produce strong oscillations local to the contact points, indicating that a reconstruction in depth is necessary local to these regions. Including the Froude number considerations by \cref{eqn:B_fastfluid} produces simulations that converge at $\order{J^{-3/2}}$ in wet regions, as shown in \cref{fig:TP_lakeoscillate_convergence}, and at $\order{J^{-2}}$ in dry regions. 

Comparing the different reconstructions (\cref{fig:TP_lakeoscillate_comparison}), we see that SkK and Bo13 both produce accurate results, while Ku07 leaves a layer of fluid on the slope and SkK diffuses a layer of fluid forward. However, the simulation Bo13 failed at time $t=2.4$ (2 \sf), at which the fluid reverses direction, due to the large value of $q_j$ that remains across the region $0<x<2$ despite fluid no longer being there, which persists due to the desingularization \cref{eqn:desing_Boller}. it is plausible that using \cref{eqn:desing_Chertock} would permit the simulation to continue, but this is likely to increase the error up to that of Ch15. The results for Ku02 is not plotted as this simulation fails at $t=0.52$ (2 \sf) for the same reasons as \cref{tp:TP_lakerest}.
\end{testproblem}

\begin{figure}[tp!]
	\begin{center}
		\includegraphics{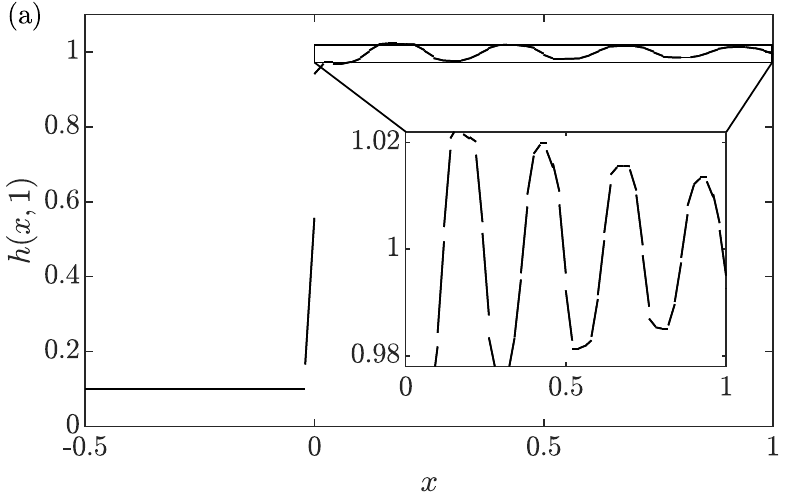}
		\includegraphics{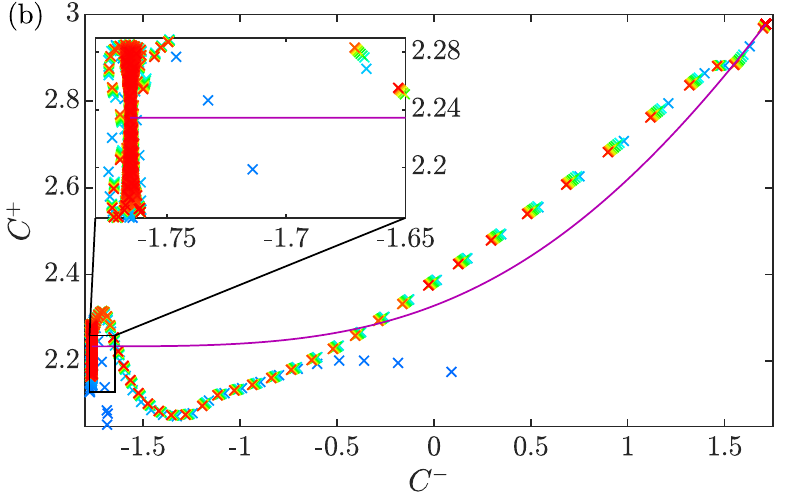}
		\includegraphics{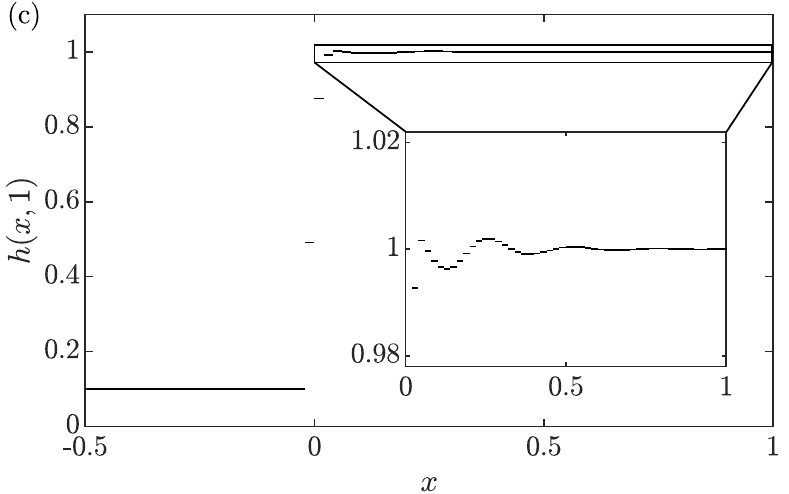}
		\includegraphics{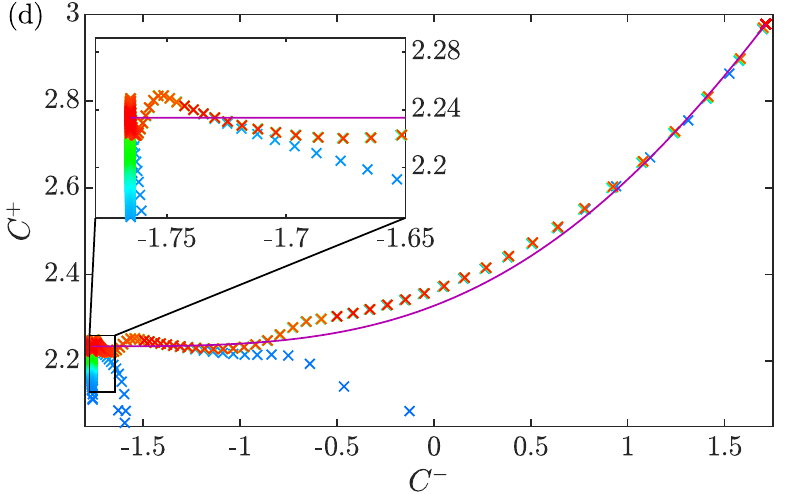}
		\includegraphics{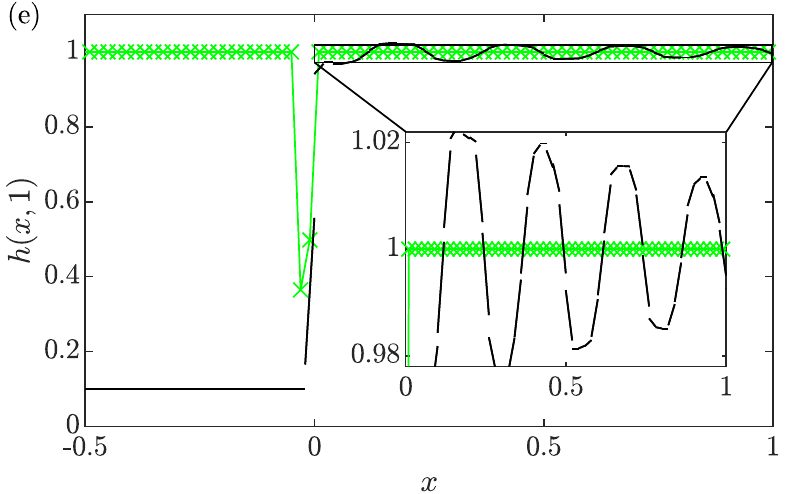}
		\includegraphics{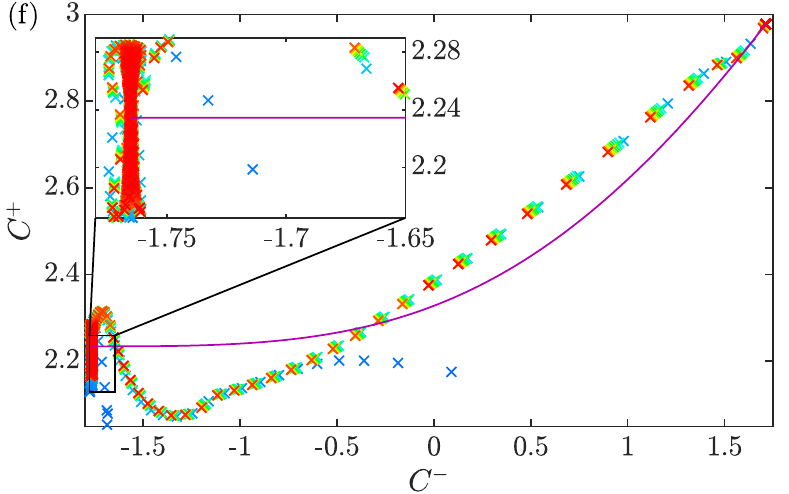}
		\includegraphics{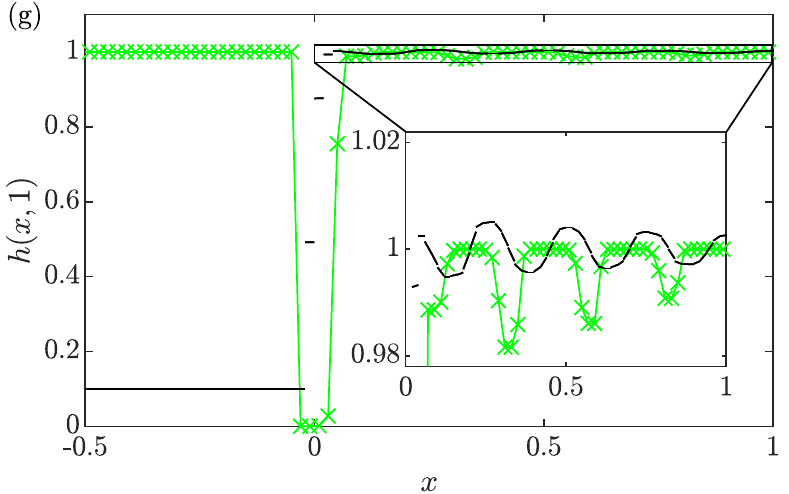}
		\includegraphics{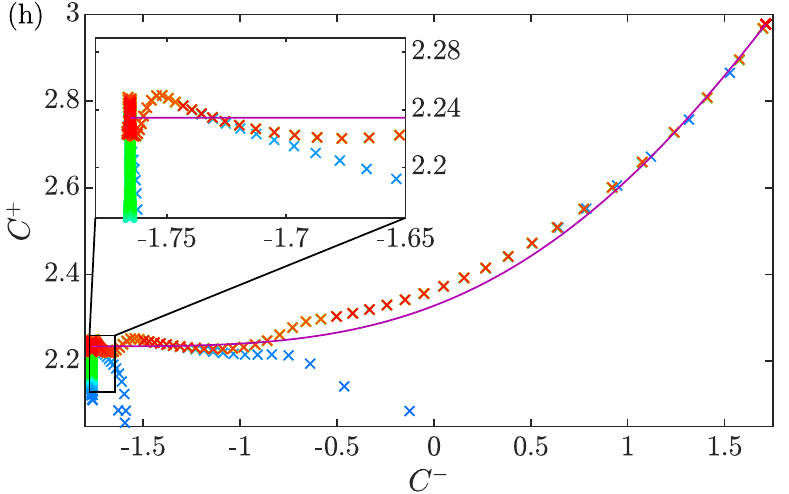}
	\end{center}
	\caption{Plots for \cref{tp:TP_slowshock}: (a,b) were simulated using $\Theta_j^h=\Theta_j^q = 1$ everywhere, that is a standard piecewise linear reconstruction without suppression; (c,d) were simulated using $\Theta_j^h=\Theta_j^q = 0$ everywhere, that is a standard piecewise constant reconstruction; (e,f) were simulated using SkK; (g,h) were simulated using SkT. Plots (a,c,e,g) show the reconstructed depth at time $t=1$, exhibiting a wave train downstream of the shock whose amplitude decays with distance. Plot (e) additionally includes the values of $\Theta_j$ and (g) the value of $\kappa_j$ for each cell as green crosses. Plots (b,d,f,h) show the cell values in every cell at times $t \in \{-1,-0.99,-0.98,\ldots, 1\}$, with early times in blue transitioning through green and yellow to late times in red, along with the Hugoniot locus of states that can be reached from $h=1$, $u=0.2345$ by a shock (purple curve).}
	\label{fig:TP_slowshock}
\end{figure}

\begin{testproblem}[Slow shock]	\label{tp:TP_slowshock}
As discussed in \cite{ar_Arora_1997,ar_Karni_1997}, when shocks move very slowly with respect to the grid cells this can cause problems in accuracy, generating a wave train. Commonly this occurs when the numerical domain is stationary ($x_L$ and $x_R$ constant), and the shock is moving at a speed with magnitude much less than $\frac*{(x_R-x_L)}{J}$, \ie taking many time steps to cross each grid cell. This is equivalent, by a change of reference frame, to a stationary shock at $x=0$ with a moving numerical domain, $\abs{\dv*{x_L}{t}}$ and $\abs{\dv*{x_R}{t}}$ being much less than $\frac*{(x_R-x_L)}{J}$. We consider such a shock, with
\begin{align}
	h &=
	\begin{cases}
		1/10	& \textrm{if } x<0,	\\
		1		& \textrm{if } x>0,
	\end{cases}
	&
	u &=
	\begin{cases}
		2.3452	& \textrm{if } x<0,	\\
		0.2345	& \textrm{if } x>0,
	\end{cases}
\end{align}
where the values of $h$ were chosen, and the corresponding values of $u$ were computed from the Rankine-Hugoniot jump conditions. Equivalently it can be written in terms of the characteristic invariants $C^+ = u+2\sqrt{h}$ and $C^-=u-2\sqrt{h}$ as
\begin{align}
	C^+ &=
	\begin{cases}
		2.9777	& \textrm{if } x<0,	\\
		2.2345	& \textrm{if } x>0,
	\end{cases}
	&
	C^- &=
	\begin{cases}
		1.7128	& \textrm{if } x<0,	\\
		-1.7655	& \textrm{if } x>0.
	\end{cases}
\end{align}
To simulate we use a grid of $J=1000$ cells on a domain of $-10+t/10 \leq x \leq 10+t/10$ over a time interval of $-1 \leq t \leq 1$. The grid cells move to the right, which is equivalent to the shock moving to the left, at a speed of $5$ grid cells per unit time. We compare different reconstruction approaches.

Firstly, we simulate using the piecewise linear $\minmod$ reconstruction without suppression (\ie we set $\Theta_j^h = \Theta_j^q = 1$). As shown in \cref{fig:TP_slowshock}(a), this yields a simulation with a relatively high amplitude wave train that propagates downstream with a large decay length due to the low diffusivity of the numerical scheme. \Cref{fig:TP_slowshock}(b) shows the evolution of the entire solution in $C^-$-$C^+$ space, using the same approach as \cite{ar_Arora_1997}. We see that the transition from the left to right state that the cells undergo follows a curve significantly perturbed from the analytic Hugoniot locus (we expect that for accurate simulation this curve should be tracked closely), and that the oscillations about the right state that are produced are almost entirely confined to the $C^+$ field (these are the oscillations we aim to reduce). 

These oscillations are significantly reduced by using a piecewise constant reconstruction (\ie we set $\Theta_j^h = \Theta_j^q = 0$). As shown in \cref{fig:TP_slowshock}(c,d), the oscillations are of a greatly decreased magnitude and decay much more rapidly. The simulation also stays much closer to the Hugoniot locus throughout the transition between the states. While using a piecewise constant reconstruction works well for simulating a single isolated shock, it will cause the scheme in smooth areas to be low accuracy.

The use of the SkK reconstruction (designed simply to permit bounds on velocity) does not have a significant effect on the wave train, \cref{fig:TP_slowshock}(e,f). This can be improved by reducing $K_{j+1/2}^\pm$, but already the reconstructed gradients are suppressed in regions where $\dv*{h}{x} \leq 10 h_j / (x_R-x_L)$ and reducing this further will excessively increase the error in regions with steep gradients.

The purpose of $\Theta_j$ as defined in \cref{eqn:SS_detector_field} is to transition between piecewise linear and piecewise constant reconstructions so that we use the piecewise liner reconstruction in regions where the solution is continuous (or at least, has no shock), and if it is differentiable the scheme has second order accuracy \cite{ar_Kurganov_2001}, whilst using piecewise constant reconstruction around shocks to suppress oscillations. The result is plotted in \cref{fig:TP_slowshock}(g,h). We see that the amplitude of the oscillations is similar in magnitude to the results in \cref{fig:TP_slowshock}(c). Indeed, both have an exponential envelope around the oscillations of the form $1 \pm 0.006 \exp(-x/\chi)$, but the piecewise linear reconstruction has a $\chi$ much larger than for the piecewise constant reconstruction due to the reduced diffusion. Additionally, the transition between the states deviates from the Hugoniot locus no more than for the piecewise constant reconstruction, in fact \cref{fig:TP_slowshock}(d) and (f) are almost indistinguishable. We conclude then that our oscillation suppression has been successful.
\end{testproblem}

\begin{figure}
	\centering
	\includegraphics{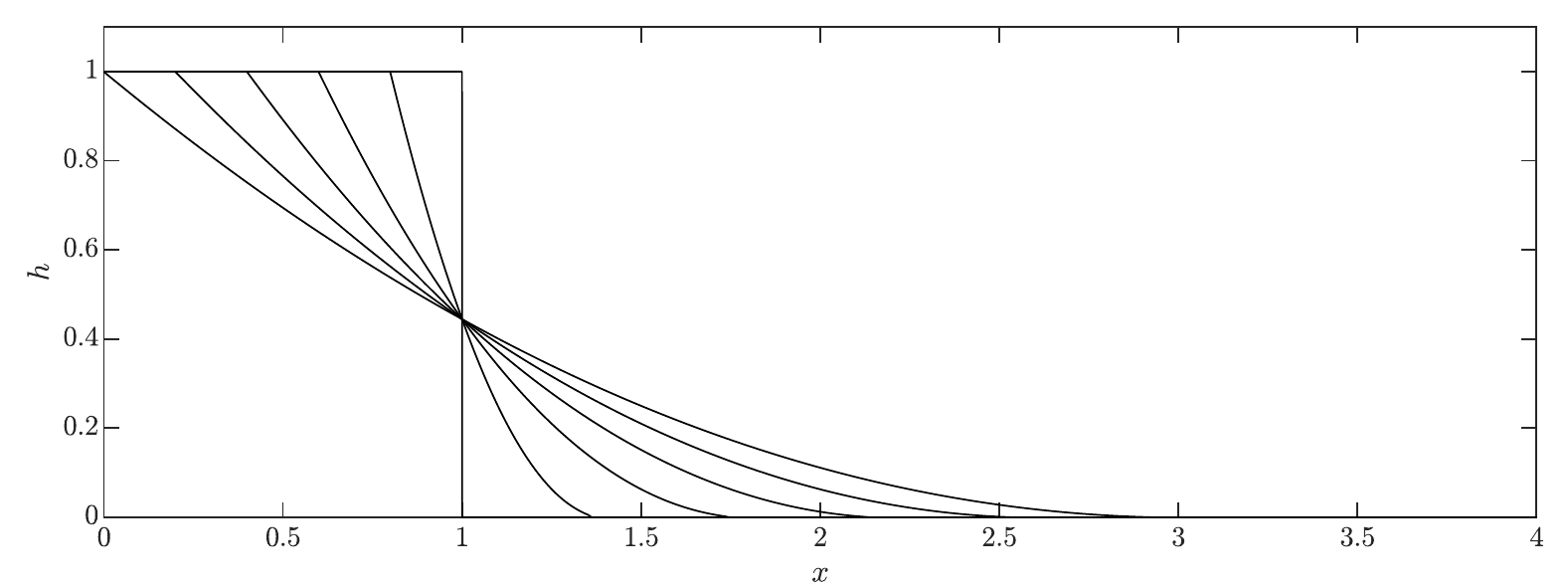}
	\caption{The depth field as a function of $x$ for \cref{tp:TP_dambreak} at times $t \in \{ 0, 0.2, \ldots , 1 \}$, taken from a simulation using SkT at resolution $J = 10^4$.}
	\label{fig:TP_dambreak}
\end{figure}
\begin{figure}
	\centering
	\includegraphics{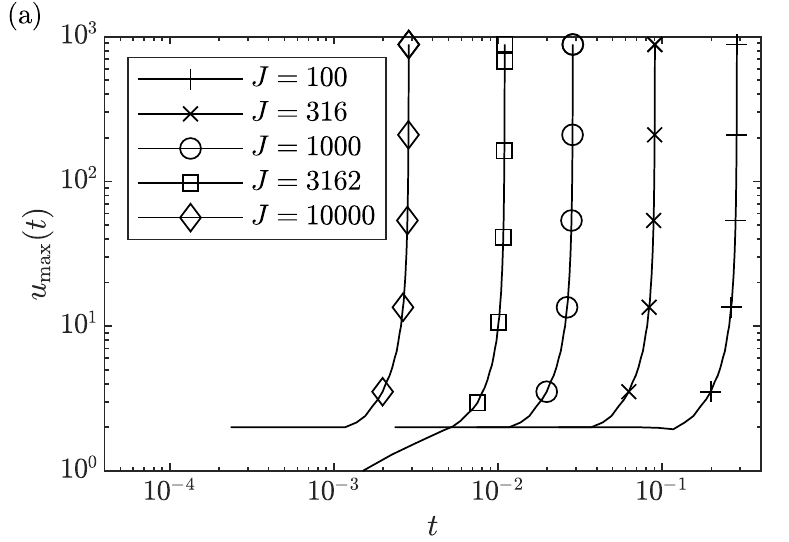}
	\includegraphics{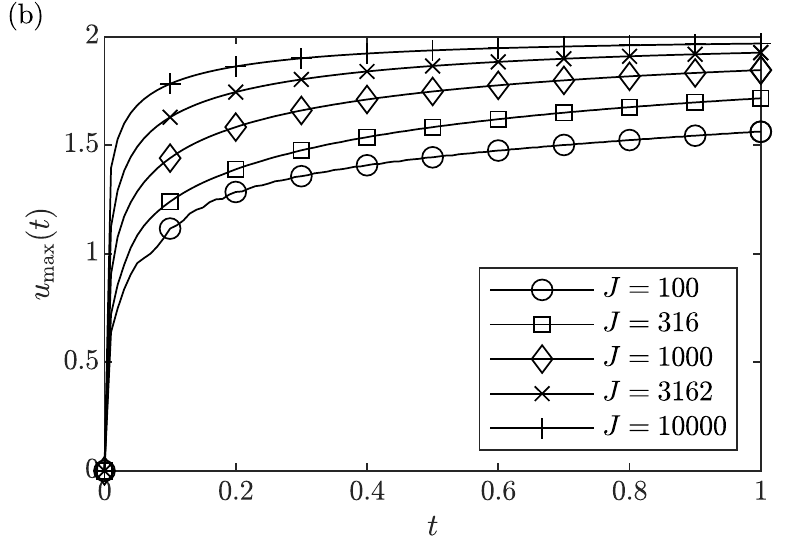}
	\caption{The maximal velocity of the fluid $u_{\max}(t) = \max_x u(x,t)$ as a function of time for a selection of resolutions. For (a) the simulations were run with $\Theta_j=1$ everywhere, and halted when $u_{\max} >10^3$. For (b) the simulations were run with SkT.}
	\label{fig:TP_dambreak_umax}
\end{figure}
\begin{figure}
	\centering
	\includegraphics{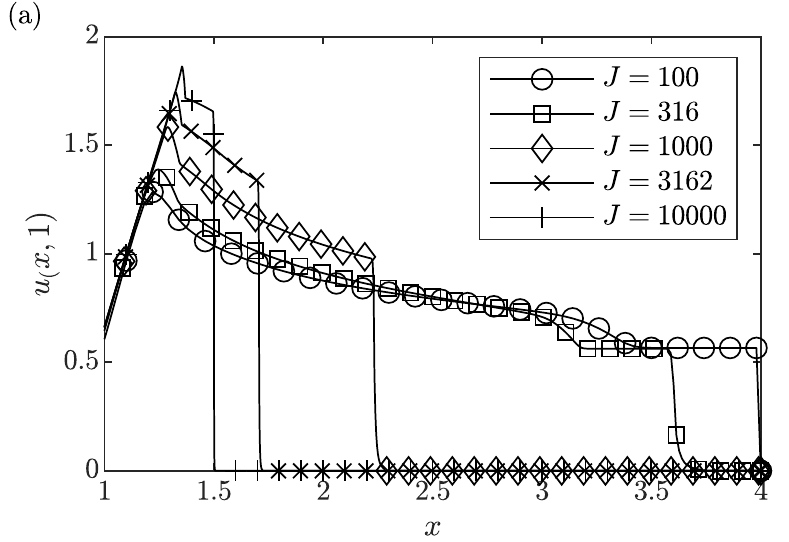}
	\includegraphics{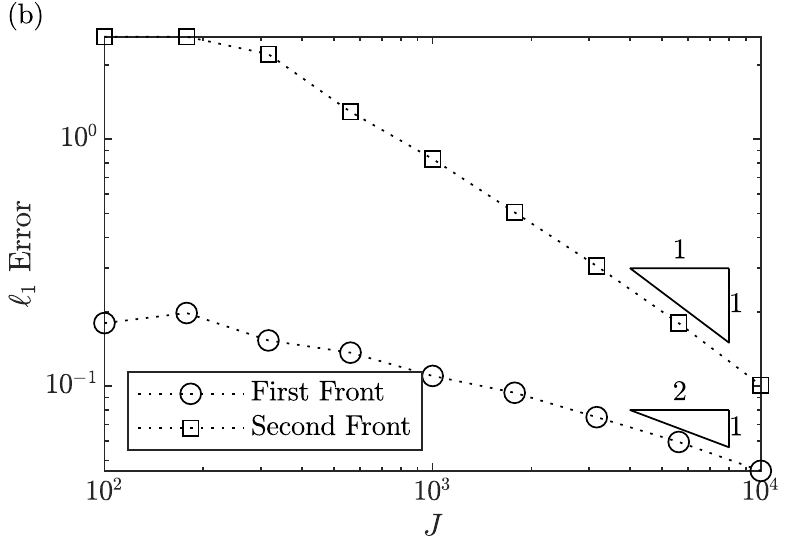}
	\includegraphics{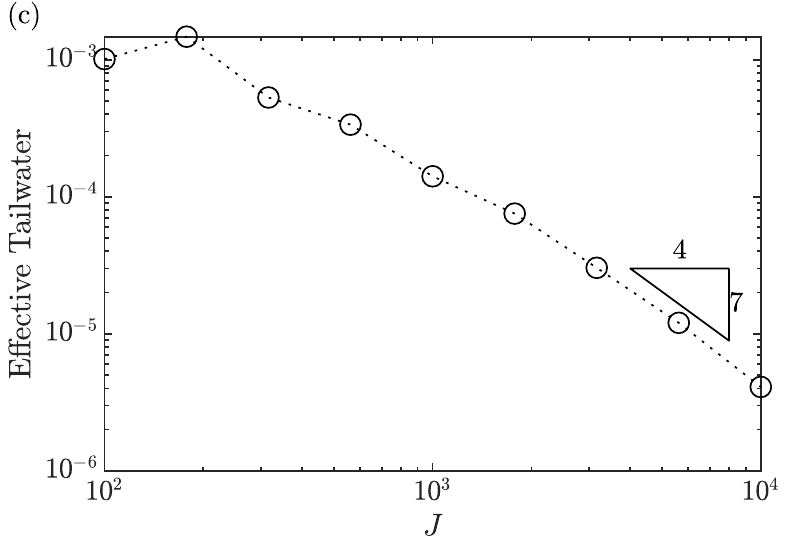}
	\caption{In (a) the velocity field at $t=0.2$ for simulations using SkT at spatial resolutions $J$ are plotted, demonstrating the two locations of the numerical front, one on either side of the location in the exact solution. Plot (b) shows error in the location of the first ($x<1.4$) and second ($x>1.4$) numerical fronts, compared to the exact front location ($x=1.4$). In (c) we plot the tailwater depth that would be required to cause deviation in the exact solution at the same location as the first numerical front.}
	\label{fig:TP_dambreak_front}
\end{figure}
\begin{figure}
	\centering
	\includegraphics{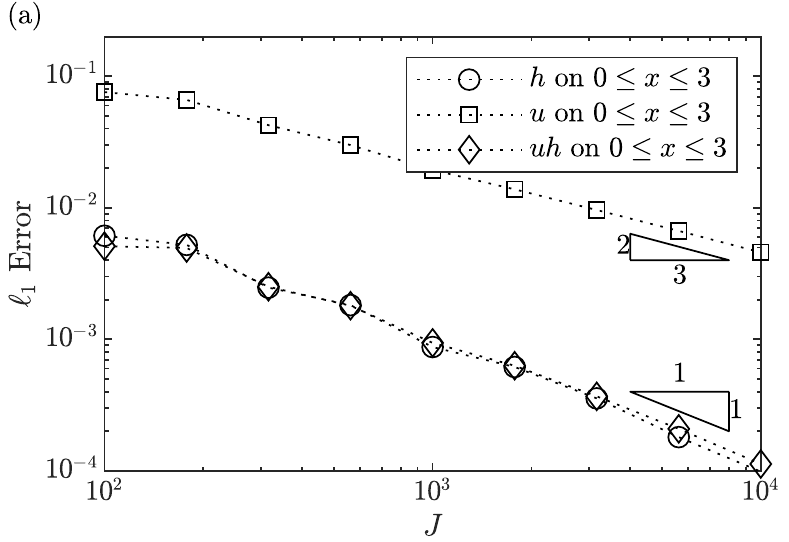}
	\includegraphics{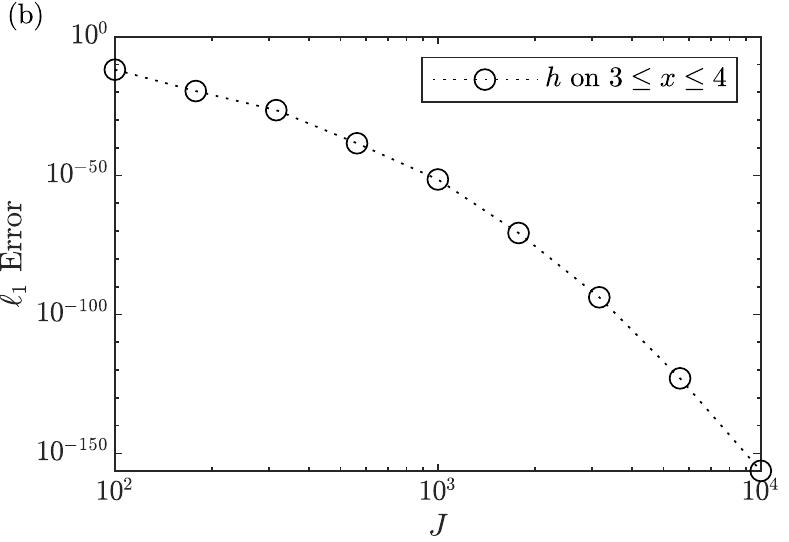}
	\caption{The $\ell_1$ error between the simulations using SkT and the exact solution as a function of the resolution $J$ for \cref{tp:TP_dambreak}. In (a) the error is evaluated at time $t=1$ for the wetted region $0 \leq x \leq 3$, and in (b) for the dry region $3 \leq x \leq 4$.}
	\label{fig:TP_dambreak_convergence}
\end{figure}

\begin{testproblem}[Dam break]	\label{tp:TP_dambreak}
As a final test problem we consider the standard dam break problem on $0 \leq x \leq 4$, which has initial conditions $u(x,0)=0$ everywhere and $h(x,0)=1$ for $x \leq 1$ and $h(x,0)=0$ for $x>1$, and we impose $u=0$ at $x \in \{0,4\}$. The exact solution to this problem for $t \leq 1$ is 
\begin{align}	\label{eqn:TP_dambreak_soln}
	h(x,t) &= 
	\begin{cases}
		1 											& \text{for } x-1 \leq -t,			\\
		\ppar*{\dfrac{2}{3} - \dfrac{x-1}{3t}}^2	& \text{for } -t \leq x-1 \leq 2t,	\\
		0											& \text{for } 2t \leq x-1,
	\end{cases}
	&
	u(x,t) &= 
	\begin{cases}
		0 											& \text{for } x-1 \leq -t,			\\
		\dfrac{2}{3} + \dfrac{2(x-1)}{3t}			& \text{for } -t \leq x-1 \leq 2t,	\\
		\text{undefined}							& \text{for } 2t \leq x-1.
	\end{cases}
\end{align}
See \cref{fig:TP_dambreak} for a plot of the depth field. Attempting to simulate from this initial condition using the piecewise linear minmod reconstruction without suppression (\ie setting $\Theta_j = 1$) causes singular velocities, as shown in \cref{fig:TP_dambreak_umax}(a). This is in part due to the lack of an upper bound for velocity for cells much deeper than their neighbour(s), see \cite{ar_Skevington_N002_Recon_Theory}, but also the effect of numerical diffusion transporting too much inertia to the front. The SkK scheme actually makes this worse, as it increases the diffusion of inertia without increasing the diffusion of depth, but if $\Theta_j^h = \Theta_j^q = \kappa_j$ is used then the velocities should remain bounded. The absence of singular velocities when using SkK is demonstrated clearly by \cref{fig:TP_dambreak_umax}(b), the velocities are bounded by and tending towards $2$, which is the exact solution at $x=2t+1$.

The velocity field at $t=0.2$ is plotted in \cref{fig:TP_dambreak_front}(a), at which time the exact solution has its peak velocity at the front $x = x_f = 1.4$. Numerically we find that the peak velocity occurs at some location $x_{f1}<x_f$ which we term the first front. Beyond this location, the velocity slowly decays up to some location $x_{f2}>x_f$ at which the velocity rapidly deceases to zero. We term this location the second front, and it is the greatest extent of the fluid layer. The deviation of these numerical front locations from the location of the front in the exact solution is shown in \cref{fig:TP_dambreak_front}(b), from which we see that $\abs{x_{f1}-x_f} = \order{J^{-1/2}}$ and $\abs{x_{f2}-x_f} = \order{J^{-3/2}}$.

We compare our approach to an alternative method that could be used to suppress these singular velocities, which is to take $\Theta_j = 1$ and impose a tail water. That is, simulate with an initial condition of $h=h_r \ll 1$ on $x>1$. This is then a Stoker problem \cite{bk_Stoker_WW}, see \cite{ar_Delestre_2016} for a discussion. To compare our method to the tailwater approach we compute the depth of a tailwater that would be required for the exact solution to the Stoker problem to deviate from \cref{eqn:TP_dambreak_soln} at $x = x_{f1}(t)$.  By the results in \cite{ar_Delestre_2016} the Stoker problem's solution is equal to \cref{eqn:TP_dambreak_soln} for $x \leq x_{\textrm{st}}$, where
\begin{align}
	x_{\textrm{st}} - 1 &= t (2 - 3 c_m),
	&\textrm{and $c_m$ satisfies}&&
	- 8 h_r c_m^2 (1 - c_m)^2 + (c_m^2 - h_r)^2 (c_m^2 + h_r) &= 0.
\end{align}
To construct an approximation to the solution we set $c_m = \epsilon$, $h_r \sim H \epsilon^p$. The requirement $c_m^2 > h_r$ is required by the entropy conditions on the shock, thus $p \geq 2$. Balance in the implicit equation for $c_m$ gives $p=4$, $H=1/8$, thus $h_r \sim c_m^4 / 8$ and the effective tailwater depth for our scheme is
\begin{equation}
	h_r \sim \frac{1}{8} \ppar*{\frac{2}{3} - \frac{x_{f1} - 1}{3t}}^4.
\end{equation}
The computed effective tailwater depths are plotted in \cref{fig:TP_dambreak_front}(c) and show that $h_r = \order{J^{-7/4}}$. This demonstrates that our method is equivalent, in terms of accuracy, to a tailwater that tends to zero at almost second order in resolution.

Examining the convergence of the scheme, \cref{fig:TP_dambreak_convergence}, we see that $h$ and $uh$ converge at first order, which is the best that can be expected for a discontinuous initial condition. This is better than the convergence that would be seen if simulating using a tailwater, which would be $\order{J^{-1/4}}$ as per the previous discussion. In the region that is dry in the exact solution we find that the volume there converges extremely rapidly.
\end{testproblem}
\section{Summary}	\label{sec:END}

We have demonstrated a modified version of the well-balanced reconstruction for the shallow water equations presented in \cite{ar_Skevington_N002_Recon_Theory}. The reconstruction is a convex combination of reconstructions in depth and in surface elevation. This allows us to reconstruct in depth when the fluid is shallow in comparison to bed variation or at high Froude number, and in surface elevation when the fluid is deep and at low Froude number. In \cref{sec:TP} the reconstruction was shown to accurately reproduce steady states (\cref{tp:TP_lakerest}), draining (\cref{tp:TP_lakedrain}), and deep dynamic states (\cref{tp:TP_lakeoscillate}). For resolving pure steady states, \cite{ar_Bollermann_2013} was found to be the superior scheme. However, our scheme was found to be superior for time evolving situations, such as draining a slope or an oscillating lake, which appear in simulations as transient phenomena on the path to steady state.

Additionally, in \cref{sec:SS} we showed how the reconstruction could be modified to suppress oscillations in slow shocks. This modification consisted of suppressing the reconstructed gradients, potentially all the way to a piecewise constant reconstruction, so that shocks could be resolved accurately. The success of this approach was demonstrated in \cref{tp:TP_slowshock}. We find that this gradient suppression is able to assist in obtaining bounded velocities for the dam break problem, see \cref{tp:TP_dambreak}. Crucially, this approach is designed to be compatible with well balancing, and has been demonstrated as such.

We conclude that the approach to reconstruction presented here is suitable for a wide class of systems.
\section*{Acknowledgements}

This work was supported by the EPSRC [grant number EP/M506473/1]. The author would also like to thank A. J. Hogg for his constructive comments regarding drafts of this article.

\bibliographystyle{plain}
\bibliography{Bibliography/Mine,Bibliography/NumericalMethodsHyperbolic,Bibliography/Books,Bibliography/ShallowWater,Bibliography/AnalysisHyperbolic}

\end{document}